# Feedback Stabilization of Tank-Liquid System with Robustness to Wall Friction

**Iasson Karafyllis[*], Filippos Vokos[*] and Miroslav Krstic[**]**

[*]Dept. of Mathematics, National Technical University of Athens,
Zografou Campus, 15780, Athens, Greece,
emails: iasonkar@central.ntua.gr, fivo33@hotmail.com

[**]Dept. of Mechanical and Aerospace Eng., University of California, San Diego, La Jolla, CA 92093-0411, U.S.A., email: krstic@ucsd.edu

**Abstract**

We solve the feedback stabilization problem for a tank, with friction, containing a liquid modeled by the viscous Saint-Venant system of Partial Differential Equations. A spill-free exponential stabilization is achieved, with robustness to the wall friction forces. A Control Lyapunov Functional (CLF) methodology with two different Lyapunov functionals is employed. These functionals determine specific parameterized sets which approximate the state space. The feedback law is designed based only on one of the two functionals (which is the CLF) while the other functional is used for the derivation of estimates of the sup-norm of the velocity. The feedback law does not require the knowledge of the exact relation of the friction coefficient. Two main results are provided: the first deals with the special case of a velocity-independent friction coefficient, while the second deals with the general case. The obtained results are new even in the frictionless case.

## 1. Introduction

The Saint-Venant model or shallow water model is a well-known mathematical model which has been used extensively since its first derivation in [1]. Many modifications of this model take into account various types of forces like gravity, viscous stresses, surface tension, friction forces (see [5,6,7,17,24,25,27,31]).

The stabilization of fluid systems is an important and challenging problem (see for instance [32]). Control studies of the Saint-Venant model focus on the inviscid Saint-Venant model (i.e., the model that ignores viscous stresses and surface tension) and its linearization around an equilibrium point (see [2,3,4,8,9,10,11,12,13,14,15,18,19,26,29]). The feedback design is performed by employing either the backstepping methodology (see [13,14]) or the Control Lyapunov Functional (CLF) methodology (see [2,3,4,9,10,11,12,15,18,19,29]). The only works that study the feedback stabilization problem for the viscous Saint-Venant model without linearization around an equilibrium point are [22,23]. More specifically, in [22,23] a feedback control law is constructed by employing the CLF methodology for the viscous Saint-Venant model without wall friction and surface tension (state feedback in [22] and output feedback in [23]). In addition [22,23] are the only works that guarantee a spill-free movement of the fluid.

Friction forces appear inevitably as the result of the interaction between the fluid and the walls of the tank that contain the fluid. This interaction is complicated and empirical relations have been proposed for the description of the forces that are developed due to wall friction (see [30]). In the mathematical literature, researchers have used various relations in order to express the effect of wall friction; see for example [4,5,15,17,18]. It is clear that wall friction forces are always present and sometimes their effect is not negligible.

In this paper we solve the feedback stabilization problem for a tank, with friction, containing a liquid modeled by the viscous Saint-Venant system of Partial Differential Equations (PDEs). A spill-free exponential stabilization is achieved, with robustness to the wall friction forces. A CLF methodology with two different Lyapunov functionals is employed. These functionals determine specific parameterized sets which approximate the state space. The feedback law is designed based only on one of the two functionals (which is the CLF) while the other functional is used for the derivation of estimates of the sup-norm of the velocity. The feedback law does not require the knowledge of the exact relation of the friction coefficient. Two main results are provided: the first (Theorem 1) deals with the special case of a velocity-independent friction coefficient, while the second (Theorem 2) deals with the general case. The obtained results are new even in the frictionless case (Corollary 1).

The paper is structured as follows. In Section 2 we describe the control problem and its main objective. In Section 3 we present the statements of the main results of this work. We also define the state space of the control problem and we present functionals which are essential for the derivation of the stability estimates. This section also includes auxiliary results which are used for the proofs of the main results. Section 4 is devoted to the proofs of all the results that we present in Section 2. Finally, in Section 5 we give the conclusions of this work and suggest topics for future research.

**Notation.** Throughout this paper, we adopt the following notation.

* $\Re_+ = [0, +\infty)$ denotes the set of non-negative real numbers.

* Let $S \subseteq \Re^n$ be an open set and let $A \subseteq \Re^n$ be a set that satisfies $S \subseteq A \subseteq cl(S)$. By $C^0(A; \Omega)$, we denote the class of continuous functions on $A$, which take values in $\Omega \subseteq \Re^m$. By $C^k(A; \Omega)$, where $k \geq 1$ is an integer, we denote the class of functions on $A \subseteq \Re^n$, which takes values in $\Omega \subseteq \Re^m$ and has continuous derivatives of order $k$. In other words, the functions of class $C^k(A; \Omega)$ are the functions which have continuous derivatives of order $k$ in $S = \text{int}(A)$ that can be continued continuously to all points in $\partial S \cap A$. When $\Omega = \Re$ then we write $C^0(A)$ or $C^k(A)$. When $I \subseteq \Re$ is an interval and $G \in C^1(I)$ is a function of a single variable, $G'(h)$ denotes the derivative with respect to $h \in I$.

* Let $I \subseteq \Re$ be an interval, let $a < b$ be given constants and let $u : I \times [a,b] \to \Re$ be a given function. We use the notation $u[t]$ to denote the profile at certain $t \in I$, i.e., $(u[t])(x) = u(t,x)$ for all $x \in [a,b]$. When $u(t,x)$ is (twice) differentiable with respect to $x \in [a,b]$, we use the notation $u_x(t,x)$ ($u_{xx}(t,x)$) for the (second) derivative of $u$ with respect to $x \in [a,b]$, i.e., $u_x(t,x) = \frac{\partial u}{\partial x}(t,x)$ ($u_{xx}(t,x) = \frac{\partial^2 u}{\partial x^2}(t,x)$). When $u(t,x)$ is differentiable with respect to $t$, we use the notation $u_t(t,x)$ for the derivative of $u$ with respect to $t$, i.e., $u_t(t,x) = \frac{\partial u}{\partial t}(t,x)$.

* Given a set $U \subseteq \Re^n$, $\chi_U$ denotes the characteristic function of $U$, i.e. the function defined by $\chi_U(x) := 1$ for all $x \in U$ and $\chi_U(x) := 0$ for all $x \notin U$. The sign function $\text{sgn}: \Re \to \Re$ is the function defined by the relations $\text{sgn}(x) = 1$ for $x > 0$, $\text{sgn}(0) = 0$ and $\text{sgn}(x) = -1$ for $x < 0$.

* Let $a < b$ be given constants. For $p \in [1, +\infty)$, $L^p(a,b)$ is the set of equivalence classes of Lebesgue measurable functions $u: (a,b) \to \Re$ with $\|u\|_p := \left( \int_a^b |u(x)|^p \, dx \right)^{1/p} < +\infty$. $L^\infty(a,b)$ is the set of equivalence classes of Lebesgue measurable functions $u: (a,b) \to \Re$ with $\|u\|_\infty := \underset{x \in (a,b)}{\text{ess sup}}(|u(x)|) < +\infty$. For an integer $k \geq 1$, $H^k(a,b)$ denotes the Sobolev space of functions in $L^2(a,b)$ with all its weak derivatives up to order $k \geq 1$ in $L^2(a,b)$.

## 2. Description of the Problem

We consider a one-dimensional model for the motion of a tank. The tank contains a viscous, Newtonian, incompressible liquid. The presence of viscosity, denoted by $\mu$, is crucial in our approach to the problem. We rely on the viscosity for establishing a region of attraction and we use the viscosity as a gain in the controller on the difference between the boundary liquid levels. The tank is subject to a force that can be manipulated. We assume that the liquid pressure is hydrostatic and consequently, the liquid is modeled by the one-dimensional (1-D) viscous Saint-Venant equations, whereas the tank obeys Newton's second law and consequently we consider the tank acceleration to be the control input.

The control objective is to drive asymptotically the tank to a specified position without the liquid spilling out and having both the tank and the liquid within the tank brought to rest. Let the position of the left side of the tank at time $t \geq 0$ be $a(t)$ and let the length of the tank be $L > 0$ (a constant). The equations describing the motion of the liquid within the tank are

$$H_t + (Hu)_z = 0, \text{ for } t > 0, \ z \in [a(t), a(t) + L] \tag{2.1}$$

$$\begin{aligned}(Hu)_t &+ \left( Hu^2 + \frac{1}{2} gH^2 \right)_z = \mu (Hu_z)_z \\ &- \kappa\big(H(t, a(t)+x), u(t, a(t)+x) - \dot{a}(t)\big) u(t, a(t)+x), \\ &- \kappa\big(H(t, a(t)+x), u(t, a(t)+x) - \dot{a}(t)\big) \dot{a}(t) \\ & \quad \text{for } t > 0, \ z \in (a(t), a(t)+L) \end{aligned} \tag{2.2}$$

where $H(t,z) > 0$, $u(t,z) \in \Re$ are the liquid level and the liquid velocity, respectively, at time $t \geq 0$ and position $z \in [a(t), a(t)+L]$, $\kappa \in C^0((0, +\infty) \times \Re; \Re_+)$ is the friction coefficient that depends on the liquid level and the relative velocity of the fluid with respect to the tank, while $g, \mu > 0$ (constants) are the acceleration of gravity and the kinematic viscosity of the liquid, respectively.

Various (empirical) relations have been used for the friction coefficient in the literature:
* in [4,18] the authors use the relation $\kappa(h,v) = c_f |v|$, where $c_f > 0$ is a constant,
* in [5] the authors use the relation $\kappa(h,v) = r_0 + r_1 h |v|$, where $r_0 > 0$, $r_1 \geq 0$ are constants,

- in [15] the authors use the relation $\kappa(h,v) = rh^{-1/3}(b+2h)^{4/3}|v|$, where $r > 0$ is a constant and $b > 0$ is the (constant) width of the tank,
- in [17] the authors derive the velocity-independent relation $\kappa(h) = \dfrac{3\mu c}{3\mu + 4ch}$, where $c > 0$ is a constant.

It is clear that there is considerable uncertainty in the friction coefficient.

The liquid velocities at the walls of the tank must coincide with the tank velocity, i.e., we have:

$$u(t, a(t)) = u(t, a(t) + L) = w(t), \text{ for } t \geq 0 \qquad (2.3)$$

where $w(t) = \dot{a}(t)$ is the velocity of the tank at time $t \geq 0$. Moreover, since the tank acceleration is the control input, we get

$$\ddot{a}(t) = -f(t), \text{ for } t > 0 \qquad (2.4)$$

where $-f(t)$, the control input to the problem, is equal to the force exerted on the tank at time $t \geq 0$ divided by the total mass of the tank. The conditions for avoiding the liquid spilling out of the tank are:

$$\begin{aligned} H(t, a(t)) &< H_{\max} \\ H(t, a(t) + L) &< H_{\max} \end{aligned} \qquad (2.5)$$

where $H_{\max} > 0$ is the height of the tank walls. Applying the transformation

$$\begin{aligned} v(t, x) &= u(t, a(t) + x) - w(t) \\ h(t, x) &= H(t, a(t) + x) \\ \xi(t) &= a(t) - a^* \end{aligned} \qquad (2.6)$$

where $a^* \in \Re$ is the specified position (a constant) to which we want to bring (and maintain) the left side of the tank, we obtain the model:

$$\dot{\xi} = w, \ \dot{w} = -f, \text{ for } t \geq 0 \qquad (2.7)$$

$$h_t + (hv)_x = 0, \text{ for } t > 0, \ x \in [0, L] \qquad (2.8)$$

$$(hv)_t + \left(hv^2 + \frac{1}{2}gh^2\right)_x = \mu(hv_x)_x - \kappa(h, v)v + hf, \text{ for } t > 0, \ x \in (0, L) \qquad (2.9)$$

$$v(t, 0) = v(t, L) = 0, \text{ for } t \geq 0 \qquad (2.10)$$

where the control input $f$ appears additively in the second equation of (2.7) and multiplicatively in (2.9). Moreover, the conditions (2.5) for avoiding the liquid spilling out of the tank become:

$$\max(h(t,0), h(t,L)) < H_{\max}, \text{ for } t \geq 0 \text{ --- Condition for no spilling out} \qquad (2.11)$$

We consider classical solutions for the PDE-ODE system (2.7)-(2.10), i.e., we consider functions $\xi \in C^2(\Re_+)$, $w \in C^1(\Re_+)$, $h \in C^1(\Re_+ \times [0, L]; (0, +\infty)) \cap C^2((0, +\infty) \times (0, L))$, $v \in C^0(\Re_+ \times [0, L]) \cap C^1((0, +\infty) \times [0, L])$ with $v[t] \in C^2((0, L))$ for each $t > 0$, which satisfy equations (2.7)-(2.10) for a given input $f \in C^0(\Re_+)$.

Using (2.8) and (2.10), we prove that for every solution of (2.7)-(2.10) it holds that $\frac{d}{dt}\left(\int_0^L h(t,x)dx\right) = 0$ for all $t > 0$. Therefore, the total mass of the liquid $m > 0$ is constant. Therefore, without loss of generality, we assume that every solution of (2.7)-(2.10) satisfies the equation

$$\int_0^L h(t,x)dx \equiv m \qquad (2.12)$$

The open-loop system (2.7)-(2.10), (2.12), i.e., system (2.7)-(2.10), (2.12) with $f(t) \equiv 0$, allows a continuum of equilibrium points, namely the points

$$h(x) \equiv h^*, \; v(x) \equiv 0, \text{ for } x \in [0, L] \qquad (2.13)$$

$$\xi \in \Re, \; w = 0 \qquad (2.14)$$

where $h^* = m/L$. We assume that the equilibrium points satisfy the conditions for no spilling out (2.11), i.e., $h^* < H_{max}$.

Our objective is to design a feedback law of the form

$$f(t) = F\big(h[t], v[t], \xi(t), w(t)\big), \text{ for } t > 0, \qquad (2.15)$$

which achieves stabilization of the equilibrium point with $\xi = 0$. Moreover, we additionally require that the "spill-free condition" (2.11) holds for every $t \geq 0$.

The existence of a continuum of equilibrium points for the open-loop system given by (2.13), (2.14) implies that the desired equilibrium point is not asymptotically stable for the open-loop system. Thus, the described control problem is far from trivial.

## 3. Construction of the feedback Law

*3.1. The Control Lyapunov Functional (CLF)*

Let $k, q > 0$ be position error and velocity gains (to be selected). Define the set $S \subset \Re^2 \times \big(C^0([0, L])\big)^2$:

$$(\xi, w, h, v) \in S \Leftrightarrow \begin{cases} h \in C^0([0, L]; (0, +\infty)) \cap H^1(0, L) \\ v \in C^0([0, L]) \\ \int_0^L h(x)dx = m \\ (\xi, w) \in \Re^2, v(0) = v(L) = 0 \end{cases} \qquad (3.1)$$

We define the following functionals for all $(\xi, w, h, v) \in S$:

$$V(\xi, w, h, v) := \delta E(h, v) + W(h, v) + \frac{qk^2}{2}\xi^2 + \frac{q}{2}(w + k\xi)^2 \qquad (3.2)$$

$$E(h, v) := \frac{1}{2}\int_0^L h(x)v^2(x)dx + \frac{g}{2}\|h - h^*\chi_{[0,L]}\|_2^2 \qquad (3.3)$$

$$W(h,v) := \frac{1}{2}\int_0^L h^{-1}(x)\big(h(x)v(x)+\mu h'(x)\big)^2 dx + \frac{g}{2}\big\|h-h^*\chi_{[0,L]}\big\|_2^2 \qquad (3.4)$$

where $\delta > 0$ and $h^* = m/L$. Moreover, define for all $(\xi, w, h, v) \in S$ with $v \in H^1(0,L)$ the functional:

$$U(\xi, w, h, v) := V(\xi, w, h, v) + \left(\frac{1}{2}\|v'\|_2^2 + \gamma V(\xi, w, h, v)\right)\exp\big(\beta V(\xi, w, h, v)\big) \qquad (3.5)$$

where $\beta, \gamma > 0$. We notice that:

- the functional $E$ is the mechanical energy of the liquid within the tank. Indeed, notice that $E$ is the sum of the potential energy ($\frac{g}{2}\|h-h^*\chi_{[0,L]}\|_2^2$) and the kinetic energy ($\frac{1}{2}\int_0^L h(x)v^2(x)dx$) of the liquid,

- the functional $W$ is a kind of mechanical energy of the liquid within the tank and has been used extensively in the literature of isentropic, compressible liquid flow (see [21,22] and references therein),

- the functional $V$ is a linear combination of the mechanical energies $E$, $W$ and a Lyapunov function for the tank ($\frac{qk^2}{2}\xi^2 + \frac{q}{2}(w+k\xi)^2$),

- the functional $U$ is a functional that can provide bounds for the sup-norm of the fluid velocity (due to the fact that $\|v'\|_2^2 \leq 2U(\xi, w, h, v)$ for all $(\xi, w, h, v) \in S$) while the other functionals cannot provide bounds for the sup-norm of the fluid velocity.

We intend to use the functional $V(\xi, w, h, v)$ defined by (3.2) as a CLF for the system. However, the functional $V(\xi, w, h, v)$ can also be used for the derivation of useful bounds for the function $h$. This is guaranteed by the following lemma.

**Lemma 1:** *Define the increasing function $G \in C^0(\Re) \cap C^1\big((-\infty,0)\cup(0,+\infty)\big)$ by means of the formula*

$$G(h) := \begin{cases} \text{sgn}(h-h^*)\left(\dfrac{2}{3}h\sqrt{h} - 2h^*\sqrt{h} + \dfrac{4}{3}h^*\sqrt{h^*}\right) & \text{for } h > 0 \\ -\dfrac{4}{3}h^*\sqrt{h^*} + h & \text{for } h \leq 0 \end{cases} \qquad (3.6)$$

*Let $G^{-1}: \Re \to \Re$ be the inverse function of $G$ and define*

$$c := \frac{1}{\mu\sqrt{\delta g}} \qquad (3.7)$$

*Then for every $(\xi, w, h, v) \in S$, the following inequality holds*

$$p_1\big(V(\xi, w, h, v)\big) \leq h(x) \leq p_2\big(V(\xi, w, h, v)\big),$$
$$\text{for all } x \in [0, L], \qquad (3.8)$$

*where the functions* $p_i : \mathfrak{R}_+ \to \mathfrak{R}$ ( $i = 1, 2$ ) *are defined by the following formulae for all* $s \geq 0$ :

$$p_1(s) := \max\left( G^{-1}(-cs), h^* - \frac{1}{\mu}\sqrt{\frac{2m(1+\delta)}{\delta}} \sqrt{s} \right)$$

$$p_2(s) := \min\left( G^{-1}(cs), h^* + \frac{1}{\mu}\sqrt{\frac{2m(1+\delta)}{\delta}} \sqrt{s} \right)$$

(3.9)

**Remark 1:** It follows from (3.9), (3.6), (3.7) and the fact that $h^* = m/L$ that $p_1(V(\xi, w, h, v)) > 0$ when $V(\xi, w, h, v) < \mu h^* \max\left( \frac{4}{3}\sqrt{\delta g h^*}, \frac{\mu \delta}{2L(1+\delta)} \right)$. Definitions (3.9) imply that $p_2 : \mathfrak{R}_+ \to \mathfrak{R}$ is an increasing function while $p_1 : \mathfrak{R}_+ \to \mathfrak{R}$ is a decreasing function.

Notice that Lemma 1 gives tighter bounds for the liquid level than Lemma 1 in [22]. Moreover, Lemma 1 in [22] can be applied only for the case $\delta = 1$, while here Lemma 1 can be applied for all $\delta > 0$.

*3.2. The state space and useful subsets of the state space*

In order to be able to satisfy the conditions for no spilling out (2.11) we need to restrict the state space. This becomes clear from the fact that the set $S$ contains states that violate the conditions for no spilling out (2.11). Define

$$X := \left\{ (\xi, w, h, v) \in S : \max_{x \in [0, L]} (h(x)) < H_{\max} \right\}$$

(3.10)

$$R := \frac{2\mu\sqrt{\delta g h^*}}{3}\left( H_{\max} - h^* \right) \min(\zeta_1, \zeta_2)$$

(3.11)

where

$$\zeta_1 := \max\left( \sqrt{\frac{H_{\max}}{h^*}} - 2\sqrt{\frac{h^*}{H_{\max}}}, \frac{3\mu\sqrt{\delta}(H_{\max} - h^*)}{4m(1+\delta)\sqrt{g h^*}} \right), \quad \zeta_2 := \frac{h^*}{H_{\max} - h^*} \max\left( 2, \frac{3\mu\sqrt{\delta} h^*}{4m(1+\delta)\sqrt{g}} \right).$$

Notice that definition (3.11), the fact that $h^* < H_{\max}$ and Lemma 1 imply that for all $(\xi, w, h, v) \in S$ with $V(\xi, w, h, v) < R$ it holds that

$$0 < p_1(V(\xi, w, h, v)) \leq h(x) \leq p_2(V(\xi, w, h, v)) < H_{\max}, \text{ for all } x \in [0, L]$$

(3.12)

Therefore, the conditions for no spilling out (2.11) are automatically satisfied when $(\xi, w, h, v) \in S$ with $V(\xi, w, h, v) < R$.

**Remark 2:** It should be noticed that $R$ depends crucially on $H_{\max}, h^*$ and $\delta$. Obviously, if $H_{\max}$ or $H_{\max} - h^*$ are small then $R$ is low because the control problem becomes more challenging and therefore the initial conditions get more restricted by the inequality $V(\xi, w, h, v) < R$. The dependence of $R$ on $\delta$ is non-monotonic and very complicated. When $\delta \to +\infty$ then we may have $R \to +\infty$ (when $H_{\max} > 2h^*$) or $R \to \dfrac{\mu^2 \left(H_{\max} - h^*\right)^2}{2m} < +\infty$ (when $H_{\max} \leq 2h^*$). The case $H_{\max} \leq 2h^*$ is the more challenging case and in this case $R$ becomes restricted when $\delta \to +\infty$.

Notice that since Lemma 1 gives tighter bounds for the liquid level than Lemma 1 in [22], inequality (3.12) is more accurate than the corresponding inequality in [22] (which can be applied only for the case $\delta = 1$).

The state space of system (2.7)-(2.10), (2.12) is considered to be the set $X$. More specifically, we consider as state space the metric space $X \subset \Re^2 \times H^1(0,L) \times L^2(0,L)$ with metric induced by the norm of the underlying normed linear space $\Re^2 \times H^1(0,L) \times L^2(0,L)$, i.e.,

$$\|(\xi, w, h, v)\|_X = \left(\xi^2 + w^2 + \|h\|_2^2 + \|h'\|_2^2 + \|v\|_2^2\right)^{1/2} \tag{3.13}$$

However, in order to construct a robust feedback law we need to approximate the state space from its interior by using certain parameterized sets that allow us to obtain useful estimates:

1) <u>Interior approximations in terms of the state variables:</u>

$$\tilde{X}(\omega) := \left\{ (\xi, w, h, v) \in S : \omega \leq \min_{x \in [0,L]} (h(x)) \leq \max_{x \in [0,L]} (h(x)) < H_{\max} \right\}, \text{ for } \omega \in \left[0, h^*\right] \tag{3.14}$$

$$\bar{X}(\omega_1, \omega_2) := \left\{ (\xi, w, h, v) \in S : \omega_1 \leq \min_{x \in [0,L]} (h(x)) \leq \max_{x \in [0,L]} (h(x)) < H_{\max}, \|v\|_\infty \leq \omega_2 \right\},$$
$$\text{for } \omega_1 \in \left[0, h^*\right], \ \omega_2 \geq 0 \tag{3.15}$$

Clearly, $\tilde{X}(0) = \bar{X}(0, +\infty) = X$ and $\bar{X}(\omega_1, \omega_2) \subseteq \tilde{X}(\omega_1) \subseteq X$ for all $\omega_1 \in \left[0, h^*\right]$, $\omega_2 \geq 0$. The physical meaning of the sets $\bar{X}(\omega_1, \omega_2)$, $\tilde{X}(\omega_1)$ is clear: the state $h$ (fluid level) is constrained to take values in $[\omega_1, H_{\max})$ and in the case of $\bar{X}(\omega_1, \omega_2)$ the state $v$ (fluid velocity) is constrained to take values in $[-\omega_2, \omega_2]$.

2) <u>Interior approximations in terms of sublevel sets of the Lyapunov functionals:</u>

$$X_V(r) := \left\{ (\xi, w, h, v) \in S : V(\xi, w, h, v) \leq r \right\}, \text{ for } r \geq 0 \tag{3.16}$$

$$X_U(r) := \left\{ (\xi, w, h, v) \in S : v \in H^1(0, L), U(\xi, w, h, v) \leq r + \gamma r \exp(\beta r) \right\},$$
$$\text{for } r \geq 0 \tag{3.17}$$

Clearly, definition (3.5) implies that $X_U(r) \subseteq X_V(r)$ for all $r \geq 0$. Moreover, inequalities (3.12) imply that

$$X_V(r) \subseteq \tilde{X}(p_1(r)) \text{ for all } r \in [0, R) \tag{3.18}$$

See also Figure 1. The following proposition allows us to guarantee an additional inclusion.

**Proposition 1:** *For every function $\varphi \in H^1(0, L)$ with $\varphi(0) = \varphi(L) = 0$ the following inequalities hold:*

$$\|\varphi\|_\infty \leq \sqrt{\frac{L}{3}} \|\varphi'\|_2 \tag{3.19}$$

*Furthermore, if $\varphi \in H^2(0, L)$ then*

$$\pi \|\varphi'\|_2 \leq L \|\varphi''\|_2 \tag{3.20}$$

Proposition 1 and definitions (3.5), (3.15), (3.17) guarantee the inclusion

$$X_U(r) \subseteq \bar{X}\left(p_1(r), \sqrt{\frac{2L}{3}(r + \gamma r \exp(\beta r))}\right) \text{ for all } r \in [0, R) \tag{3.21}$$

See also Figure 2. The following proposition guarantees that the set $X_V(r)$ for $r > 0$ contains a neighborhood of $(0, 0, h^* \chi_{[0,L]}, 0)$ (in the topology of $X$ with metric induced by the norm $\|\ \|_X$ defined by (3.13)).

**Proposition 2:** *Let $q, k, \delta > 0$ be given. Then for every $(\xi, w, h, v) \in S$ satisfying the inequality $\|(0, w, h - h^* \chi_{[0,L]}, v)\|_X \leq \varepsilon$ for some $\varepsilon > 0$ with $\varepsilon < \frac{\min(h^*, H_{\max} - h^*)}{\sqrt{L}}$, the following inequality holds:*

$$V(\xi, w, h, v) \leq \max\left(\mu^2 \left(h^* - \varepsilon \sqrt{L}\right)^{-1}, \frac{\delta + 1}{2} g, \frac{(\delta + 2) H_{\max}}{2}, q, \frac{3qk^2}{2}\right) \|(\xi, w, h - h^* \chi_{[0,L]}, v)\|_X^2 \tag{3.22}$$

*where $\|\cdot\|_X$ is defined by (3.13).*

## *3.3. Main results*

Now we are in a position to state the first main result of the present work, which deals with a special case of the friction coefficient $\kappa$; namely, the case where the function $\kappa \in C^0((0, +\infty) \times \Re; \Re_+)$ satisfies the following assumption.

**(H)** *There exists a continuous, non-increasing function $K : (0, H_{\max}] \to \Re_+$ such that for every $\omega \in (0, h^*]$ the following inequality holds:*

$$h^{-2} \kappa(h, v) \leq K(\omega), \text{ for all } (h, v) \in [\omega, H_{\max}] \times \Re \tag{3.23}$$

Assumption (H) allows us to achieve stabilization in the set $\tilde{X}(\omega)$ defined by (3.14) for every $\omega \in (0, h^*]$. More specifically, we achieve exponential stabilization in the set $X_V(r)$ defined by (3.16) for all $r \in [0, R)$ with $p_1(r) \geq \omega$ (for which $X_V(r) \subseteq \tilde{X}(\omega)$; recall inclusion (3.18)). This is shown by the following theorem.

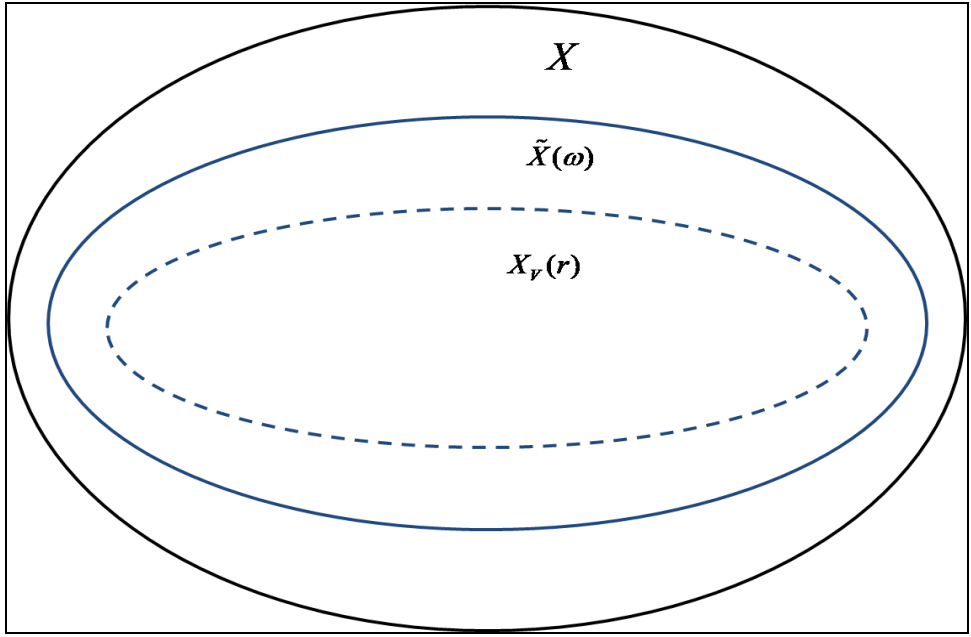

**Figure 1:** The structure of the sets $X_V(r), \tilde{X}(\omega), X$ for $r \in [0, R)$ with $p_1(r) \geq \omega$.

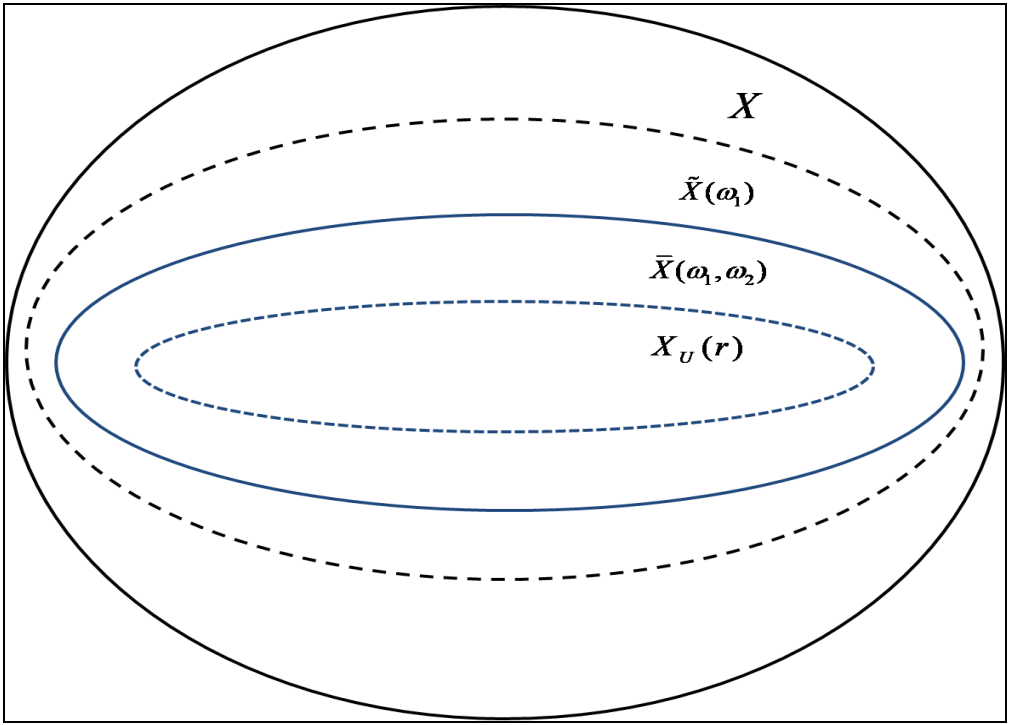

**Figure 2:** The structure of the sets $X_U(r), \bar{X}(\omega_1, \omega_2), \tilde{X}(\omega_1), X$ for $r \in [0, R)$ with $p_1(r) \geq \omega_1$ and $\sqrt{\dfrac{2L}{3}(r + \gamma r \exp(\beta r))} \leq \omega_2$.

**Theorem 1 (Special Case of Friction Coefficient-Stabilization in $X_V(r) \subseteq \tilde{X}(\omega)$):** *Suppose that the function $\kappa \in C^0((0,+\infty) \times \mathfrak{R}; \mathfrak{R}_+)$ satisfies Assumption (H). Let arbitrary $\omega \in (0, h^*]$ be given. Pick arbitrary $\delta > 0$ for which the following inequality holds:*

$$2g(\delta+1) > \mu K(\omega) \tag{3.24}$$

*Pick arbitrary $r \in [0, R)$, for which $p_1(r) \geq \omega$, where $p_1$, $R$ are defined by (3.9), (3.11). Pick arbitrary $\sigma, k, q > 0$ with*

$$k < q\theta(r) \tag{3.25}$$

*where*

$$\theta(r) := \frac{\sigma g \mu \delta \pi^2 p_1(r)}{g \mu \delta \pi^2 p_1(r) + 2\sigma L\left(mgLH_{\max}(\delta+1)^2 + 2\mu^2 \delta \pi^2 p_1(r)\right)} \tag{3.26}$$

*Then there exist constants $M, \lambda, \bar{M}, \bar{\lambda} > 0$ with the following property:*

**(P)** *Every classical solution of the PDE-ODE system (2.7)-(2.10), (2.12) and*

$$f(t) = -\sigma\left((\delta+1)\int_0^L h(t,x)v(t,x)dx + \mu(h(t,L) - h(t,0)) - q(w(t) + k\xi(t))\right), \text{ for } t > 0 \tag{3.27}$$

*with $v \in C^0(\mathfrak{R}_+; H^1(0,L)) \cap C^1((0,+\infty); H^1(0,L))$, $(\xi(0), w(0), h[0], v[0]) \in X_V(r) \subseteq \tilde{X}(\omega)$, satisfies $(\xi(t), w(t), h[t], v[t]) \in X_V(r) \subseteq \tilde{X}(\omega)$ and the following estimates for $t \geq 0$:*

$$\begin{aligned}
&\left\|(\xi(t), w(t), h[t] - h^*\chi_{[0,L]}, v[t])\right\|_X \\
&\leq M \exp(-\lambda t) \left\|(\xi(0), w(0), h[0] - h^*\chi_{[0,L]}, v[0])\right\|_X
\end{aligned} \tag{3.28}$$

$$\|v_x[t]\|_2 \leq \bar{M} \exp(-\bar{\lambda} t)\left(\left\|(\xi(0), w(0), h[0] - h^*\chi_{[0,L]}, v[0])\right\|_X + \|v_x[0]\|_2\right) \tag{3.29}$$

**Remarks on Theorem 1:**
a) Theorem 1 can be applied directly to velocity-independent friction coefficients $\kappa(h,v) = \kappa(h)$, since velocity-independent friction coefficients satisfy automatically Assumption (H). For example, the velocity-independent friction coefficient $\kappa(h) = \frac{3\mu c}{3\mu + 4ch}$, where $c > 0$ is a constant, which was derived in [17] satisfies Assumption (H) with $K(\omega) = \frac{3\mu c}{\omega^2(3\mu + 4c\omega)}$. Moreover, if there exists a constant $B > 0$ such that $\kappa(h,v) \leq B$ for all $(h,v) \in (0,+\infty) \times \mathfrak{R}$, then the friction coefficient $\kappa(h,v)$ satisfies Assumption (H) with $K(\omega) = B/\omega^2$.
b) It should be noticed that estimate (3.29) in conjunction with the inequality $\|v[t]\|_\infty \leq \sqrt{\frac{L}{3}}\|v_x[t]\|_2$ (recall Proposition 1) provides the following estimate for the sup-norm of the fluid velocity:

$$\|v[t]\|_\infty \leq \sqrt{\frac{L}{3}}\bar{M}\exp(-\bar{\lambda}t)\left(\left\|(\xi(0), w(0), h[0] - h^*\chi_{[0,L]}, v[0])\right\|_X + \|v_x[0]\|_2\right)$$

c)      Theorem 1 guarantees robust exponential stabilization of the state by means of the nonlinear feedback law (3.27). Here, the term "robust" refers to robustness with respect to friction coefficients that satisfy Assumption (H).

d)      The control parameters are $\sigma, k, q, \delta > 0$. It should be noticed that there are only two restrictions for the control parameters: (3.24) which implies that $\delta$ must be sufficiently large to overcome friction (see (3.23), (3.24)) and (3.25) which implies that the ratio $k/q$ must be sufficiently small. Definition (3.26) shows dependence of the upper bound for the ratio $k/q$ on $r, \delta, \mu, \sigma, g, L$; consequently, (3.25) is a fundamentally nonlinear gain restriction. The risk of spillage limits the gain on the tank position, even though the gain on the tank velocity is unrestricted. The risk of spillage forces us to be patient in reducing the distance to the target position.

e)      How large is the set $X_V(r)$ for which robust exponential stabilization is achieved? As we noticed above, Proposition 2 guarantees that the set $X_V(r)$ for $r > 0$ contains a neighborhood of $\left(0, 0, h^* \chi_{[0,L]}, 0\right)$ (in the topology of $X$ with metric induced by the norm $\|\ \|_X$ defined by (3.13)). However the size of the set $X_V(r)$ depends on $r \in [0, R)$ that satisfies $p_1(r) \geq \omega$ and on $\delta, q, k$ (recall definitions (3.2), (3.16)). The dependence of $X_V(r)$ on $q, k$ is clear: the larger $q$ (or $k$) the smaller the set $X_V(r)$. However, the dependence of $X_V(r)$ on $\omega \in (0, h^*]$ (for which $p_1(r) \geq \omega$ and (3.24) must hold; thus $\omega$ affects both $r$ and $\delta$), on $\delta$ (through (3.2) and through the dependence of $R$ on $\delta$; see Remark 2) is not clear: it is a very complicated, non-monotonic dependence which cannot be described easily.

f)      It should be noticed that the feedback law (3.27) does not require knowledge of the exact relation that determines the friction coefficient. The feedback law simply requires to know only a velocity-independent upper bound of the friction coefficient (see also Remark (a)). Moreover, the feedback law (3.27) does not require the measurement of the whole liquid level and liquid velocity profile and requires the measurement of only four quantities:

- the tank position $\xi(t)$ and the tank velocity $w(t)$,
- the total liquid momentum $\int_0^L h(t,x)v(t,x)dx$, and
- the liquid level difference at the tank walls $h(t,L) - h(t,0)$.

Theorem 1 shows that robustness is achieved by increasing the gain $\sigma(\delta + 1)$ of the total liquid momentum in the feedback law. However, Theorem 1 can also be applied to the frictionless case and provides a novel result even in this case.

**Corollary 1 (Frictionless Case):** *Suppose that $\kappa \equiv 0$. Pick arbitrary $\delta > 0$ and pick arbitrary $r \in [0, R)$, where $R > 0$ is defined by (3.11). Pick arbitrary $\sigma, k, q > 0$ for which (3.25) holds. Then there exist constants $M, \lambda, \bar{M}, \bar{\lambda} > 0$ for which property (P) holds.*

Corollary 1 is an improved result compared to Theorem 1 in [22] due to the following facts:
- Corollary 1 provides a larger family of feedback laws that stabilize the PDE-ODE system (2.7)-(2.10), (2.12) compared to Theorem 1 in [22], since the feedback laws provided by Theorem 1 in [22] correspond to the feedback law (3.27) with $\delta = 1$,
- Corollary 1 provides the exponential decay estimates (3.28), (3.29) while Theorem 1 in [22] provides only the exponential decay estimate (3.28).

On the other hand, Theorem 1 applies to more regular solutions than Theorem 1 in [22]. More specifically, property (P) holds for classical solutions of system (2.7)-(2.10), (2.12) and (3.27) with $v \in C^0(\Re_+; H^1(0,L)) \cap C^1((0,+\infty); H^1(0,L))$ while the requirement that $v \in C^0(\Re_+; H^1(0,L)) \cap C^1((0,+\infty); H^1(0,L))$ is absent in Theorem 1 in [22].

Our second main result of the present work deals with the general case of the friction coefficient $\kappa \in C^0((0,+\infty) \times \Re; \Re_+)$ without any further assumption. We achieve stabilization in the set $\bar{X}(\omega_1, \omega_2) \subseteq X$ defined by (3.15) for every $\omega_1 \in [0, h^*]$, $\omega_2 \geq 0$. More specifically, we achieve exponential stabilization in the set $X_U(r)$ defined by (3.17) for all $r \in [0,R)$ with $p_1(r) > \omega_1$, $\sqrt{\frac{2L}{3}}(r + \gamma r \exp(\beta r)) < \omega_2$ (for which $X_U(r) \subseteq \bar{X}(\omega_1, \omega_2)$; recall inclusion (3.21)). This is shown by the following theorem.

**Theorem 2 (General Case-Stabilization in $X_U(r) \subseteq \bar{X}(\omega_1, \omega_2)$):** *Let constants $\omega_1 \in (0, h^*)$, $\omega_2 > 0$ be given. Define*

$$\tilde{K} := \max\left\{ h^{-2}\kappa(h,v) : \omega_1 \leq h \leq H_{\max}, |v| \leq \omega_2 \right\} \tag{3.30}$$

*Pick arbitrary $\delta > 0$ that satisfies*

$$2g(\delta+1) > \mu \tilde{K} \tag{3.31}$$

*Pick arbitrary constants $\sigma, q, k > 0$ with*

$$k < q\tilde{\theta} \tag{3.32}$$

*where*

$$\tilde{\theta} := \frac{\sigma g \mu \delta \pi^2 \omega_1}{g\mu\delta\pi^2\omega_1 + 2\sigma L\left(mgLH_{\max}(\delta+1)^2 + 2\mu^2\delta\pi^2\omega_1\right)} \tag{3.33}$$

*Then for every $r \in [0, R)$ with*

$$p_1(r) > \omega_1, \quad \sqrt{\frac{2L}{3}}(r + \gamma r \exp(\beta r)) < \omega_2 \tag{3.34}$$

*where $p_1, R$ are defined by (3.9), (3.11), there exist constants $M, \lambda, \bar{M}, \bar{\lambda} > 0$ with the following property:*

**(P')** *Every classical solution of the PDE-ODE system (2.7)-(2.10), (2.12) and (3.27) with $(\xi(0), w(0), h[0], v[0]) \in X_U(r) \subseteq \bar{X}(\omega_1, \omega_2)$, satisfies $(\xi(t), w(t), h[t], v[t]) \in X_U(r) \subseteq \bar{X}(\omega_1, \omega_2)$ and estimates (3.28), (3.29) for all $t \geq 0$.*

*where $U$ is defined by (3.5) with*

$$\gamma > \frac{5(H_{\max}\tilde{K}^2 + \varepsilon_1)}{\delta\mu\tilde{\alpha}} \tag{3.35}$$

$$\beta > \max\left(\frac{4\varepsilon_2}{(2\tilde{\alpha} + \mu\delta\gamma\omega_1)\omega_1^2}, \frac{20L}{3\mu^2\delta\omega_1}\right) \tag{3.36}$$

$$\varepsilon_1 := \frac{(\delta+1)g^2}{\mu^2} H_{\max} + 3\sigma^2 L\big((\delta+1)(\delta+2)m + \delta q\big) \tag{3.37}$$

$$\varepsilon_2 := \frac{100(\delta+1)^2 R}{\delta^2 \mu^3} \tag{3.38}$$

$$\tilde{\alpha} := \frac{\min\big(\mu g, 4qk^3, 4q(q\tilde{\theta}-k), \mu\delta\big)}{2\max\big(\pi^{-2}L^2(\delta+2)\omega_1^{-1}H_{\max}, (\delta+1)gL^2 + 2\omega_1^{-1}\mu^2, qk^2, q\big)} \tag{3.39}$$

**Remarks on Theorem 2:**
a) Theorem 2 guarantees robust exponential stabilization of the state by means of the nonlinear feedback law (3.27). Again, the term "robust" refers to robustness with respect to the friction coefficient.
b) The control parameters are $\sigma, k, q, \delta > 0$. As in Theorem 1, here there are also two restrictions for the control parameters: (3.31) which implies that $\delta$ must be sufficiently large to overcome friction (see (3.30), (3.31)) and (3.32) which implies that the ratio $k/q$ must be sufficiently small. Definition (3.33) shows dependence of the upper bound for the ratio $k/q$ on $r, \delta, \mu, \sigma, g, L$; consequently, (3.32) is a fundamentally nonlinear gain restriction. Again, as we remarked in Theorem 1, the risk of spillage limits the gain on the tank position, even though the gain on the tank velocity is unrestricted.
c) How large is the set $X_U(r)$ for which robust exponential stabilization is achieved? Proposition 2 and definition (3.5) guarantees that the set $X_U(r)$ for $r > 0$ contains a neighborhood of $(0, 0, h^* \chi_{[0,L]}, 0)$ (in the topology of $\hat{X} = \{(\xi, w, h, v) \in X : v \in H^1(0, L)\}$ with metric induced by the norm $\|(\xi, w, h, v)\| = \big(\|(\xi, w, h, v)\|_X^2 + \|v'\|_2^2\big)^{1/2}$, where $\|\ \|_X$ is defined by (3.13)). The size of the set $X_U(r)$ depends on $r \in [0, R)$ that satisfies (3.34) and on $\beta, \gamma, \delta, q, k$ (recall definitions (3.2), (3.5), (3.17)). The dependence of $X_U(r)$ on $q, k$ is clear: the larger $q$ (or $k$) the smaller the set $X_U(r)$. However, the dependence of $X_U(r)$ on $\omega_1 \in (0, h^*)$ and $\omega_2 > 0$ (through (3.30), (3.31), (3.34) and (3.35)), on $\beta, \gamma, \delta$ (through (3.2), (3.5), (3.17) and through the dependence of $R$ on $\delta$; see Remark 2) is not clear: it is a very complicated, non-monotonic dependence which cannot be described easily.

In the proof of both Theorem 1 and Theorem 2, we are working with two Lyapunov functionals: $V(\xi, w, h, v)$ defined by (3.2) and $U(\xi, w, h, v)$ defined by (3.5). However, the feedback law is designed by using $V(\xi, w, h, v)$ only (it is a feedback law of $L_g V$-type). In other words, the CLF is $V(\xi, w, h, v)$ and the functional $U(\xi, w, h, v)$ is used only for the derivation of estimate (3.29) and the derivation of estimates of the sup-norm of the velocity (used in the proof of Theorem 2). Therefore, a question arises:

> "Why don't we use $U(\xi, w, h, v)$ as a CLF and design the feedback law based on $U(\xi, w, h, v)$?"

The answer is clear: if we design the feedback law (of $L_g V$-type) with $U(\xi, w, h, v)$ as a CLF then the feedback law will contain terms with spatial derivatives of the velocity. Therefore, the measurement requirements of the feedback law will be far more demanding. We avoid the measurement of spatial derivatives of the velocity by using $V(\xi, w, h, v)$ as a CLF. However, since $V(\xi,w,h,v)$ cannot give us the useful estimate (3.29) we also use $U(\xi,w,h,v)$. This situation is something that cannot happen in finite-dimensional systems: a CLF for a finite-dimensional system can provide all estimates for the solution. *Therefore, in the infinite-dimensional case we may need one functional to play the role of a CLF and a different functional for the derivation of useful stability estimates.* This situation has been observed again in [20] in the context of a nonlinear reaction-diffusion PDE, where the CLF could only provide stability estimates for the $L^2$ norm of the state and different methodologies (Stampacchia's method and an additional functional) were used for the derivation of stability estimates in the sup-norm and the $H^1$ norm of the state. The advice for the researcher that arises from these studies is clear:

- Use as a CLF a functional that can provide bounds for the state norm in a functional space of minimal regularity. Such a CLF will minimize the measurement requirements of the feedback law.
- Use different functionals in order to obtain bounds for the state norm in more demanding functional spaces.

*3.4. Auxiliary results*

For the proof of Theorem 1, we need some auxiliary lemmas. The first auxiliary lemma provides formulas for the time derivatives of the energy functionals defined by (3.3), (3.4).

**Lemma 2:** *For every classical solution of the PDE-ODE system (2.7)-(2.10), (2.12) the following equations hold for all $t > 0$:*

$$\frac{d}{dt}E(h[t],v[t]) = -\mu\int_0^L h(t,x)v_x^2(t,x)dx + f(t)\int_0^L h(t,x)v(t,x)dx$$
$$-\int_0^L \kappa(h(t,x),v(t,x))v^2(t,x)dx \quad (3.40)$$

$$\frac{d}{dt}W(h[t],v[t]) = -\mu g\|h_x[t]\|_2^2 + f(t)\int_0^L \left(h(t,x)v(t,x)+\mu h_x(t,x)\right)dx$$
$$-\int_0^L \kappa(h(t,x),v(t,x))v^2(t,x)dx - \mu\int_0^L h^{-1}(t,x)h_x(t,x)\kappa(h(t,x),v(t,x))v(t,x)dx \quad (3.41)$$

*where $E,W$ are the functionals defined by (3.3), (3.4), respectively.*

The two next auxiliary lemmas provide useful inequalities for the CLF $V$ defined by (3.2).

**Lemma 3:** *Let $q,k,\delta > 0$ be given. Then there exists a non-decreasing function $\Lambda:[0,R)\to(0,+\infty)$, where $R > 0$ is defined by (3.11) such that for every $(\xi,w,h,v)\in X$ with $v\in H^1(0,L)$ and $V(\xi,w,h,v) < R$, the following inequality holds:*

$$\frac{V(\xi,w,h,v)}{\Lambda(V(\xi,w,h,v))} \le \|h'\|_2^2 + \int_0^L h(x)(v'(x))^2 dx + \xi^2 + (w+k\xi)^2 \quad (3.42)$$

**Lemma 4:** *Let $q,k,\delta > 0$ be given. Then there exist non-decreasing functions $G_i:[0,R)\to(0,+\infty)$, $i=1,2$, where $R > 0$ is defined by (3.11), such that for every $(\xi,w,h,v)\in X$ with $V(\xi,w,h,v) < R$, the following inequality holds:*

$$\frac{V(\xi,w,h,v)}{G_2(V(\xi,w,h,v))} \leq \left\|(\xi,w,h-h^*\chi_{[0,L]},v)\right\|_X^2 \leq V(\xi,w,h,v)G_1(V(\xi,w,h,v)) \quad (3.43)$$

where $\|\cdot\|_X$ is defined by (3.13).

Inequality (3.42) provides an estimate of the dissipation rate of the Lyapunov functional for the closed-loop system (2.7)-(2.10), (2.12) with (3.27). On the other hand, inequalities (3.43) provide estimates of the Lyapunov functional in terms of the norm of the state space.

The following two lemmas are technical lemmas which provide useful differential inequalities.

**Lemma 5:** *Let $\sigma,k,q,\delta > 0$ be given constants and let $r \in [0,R)$ be a constant, where $R > 0$ is the constant defined by (3.11). Then for every classical solution of the PDE-ODE system (2.7)-(2.10), (2.12) and (3.27) with $v \in C^0(\Re_+; H^1(0,L)) \cap C^1((0,+\infty); H^1(0,L))$ the following inequalities hold for all $t > 0$ for which $V(\xi(t),w(t),h[t],v[t]) < R$:*

$$\frac{d}{dt}V(\xi(t),w(t),h[t],v[t]) \leq -\frac{\mu g}{4}\|h_x[t]\|_2^2 - q(q\theta(r)-k)(w(t)+k\xi(t))^2$$
$$-\frac{\mu\delta}{2H_{max}}\left(2H_{max} - p_1(r)\frac{p_2(V(t))}{p_1(V(t))}\right)\int_0^L h(t,x)v_x^2(t,x)dx \quad (3.44)$$
$$-\left(\delta+1-\frac{\mu\bar{K}(t)}{2g}\right)\int_0^L \kappa(h(t,x),v(t,x))v^2(t,x)dx - qk^3\xi^2(t)$$

$$\frac{d}{dt}\left(\frac{1}{2}\|v_x[t]\|_2^2\right) \leq -\frac{\mu\pi^2}{4L^2}\|v_x[t]\|_2^2 + \frac{5L}{6\mu}\|v_x[t]\|_2^4$$
$$+\frac{\varepsilon_2 V(t)}{p_1^2(V(t))}\|v_x[t]\|_2^2 + \frac{5}{\delta\mu}\left(H_{max}\bar{K}^2(t)+\varepsilon_1\right)V(t) \quad (3.45)$$

*where $V(t) = V(\xi(t),w(t),h[t],v[t])$, $\theta(r)$ is defined by (3.26), $p_i: \Re_+ \to \Re$ ($i=1,2$) are the functions defined by (3.9), $\varepsilon_1, \varepsilon_2 > 0$ are defined by (3.37), (3.38), respectively, and*

$$\bar{K}(t) := \max_{0 \leq x \leq L}\left(h^{-2}(t,x)\kappa(h(t,x),v(t,x))\right) \quad (3.46)$$

**Lemma 6:** *Let $\beta,\gamma,\sigma,k,q,\delta > 0$ be given constants. Define for all $(\xi,w,h,v) \in S$ with $v \in H^1(0,L)$ the functional $U(\xi,w,h,v)$ by means of (3.5). Let $r \in [0,R)$ be a constant, where $R > 0$ is the constant defined by (3.11) and define for $s \in [0,R)$*

$$\varphi(s) := 2H_{max}p_1(s) - p_1(r)p_2(s) \quad (3.47)$$

$$\alpha(s) := \frac{1}{\Lambda(s)}\min\left(\frac{\mu g}{4}, qk^3, q(q\theta(r)-k), \frac{\mu\delta}{4H_{max}}\left(2H_{max} - p_1(r)\frac{p_2(s)}{p_1(s)}\right)\right) \quad (3.48)$$

*where $p_i: \Re_+ \to \Re$ ($i=1,2$) are the functions defined by (3.9), $\theta(r)$ is defined by (3.26) and $\Lambda$ is the function involved in (3.42). Suppose that (3.25) holds. Then for every classical solution of the PDE-ODE system (2.7)-(2.10), (2.12) and (3.27) with $v \in C^0(\Re_+; H^1(0,L)) \cap C^1((0,+\infty); H^1(0,L))$ the* following inequality holds for all $t > 0$ for which $V(\xi(t), w(t), h[t], v[t]) < R$ and $\varphi(V(\xi(t), w(t), h[t], v[t])) > 0$:

$$\frac{d}{dt}U(\xi(t), w(t), h[t], v[t]) \leq -\alpha(V(t))V(t)$$
$$-\left(\frac{\delta\gamma}{H_{max}}\varphi(V(t)) + \frac{\pi^2}{L^2}\right)\frac{\mu}{4}\exp(\beta V(t))\|v_x[t]\|_2^2$$
$$-\left(\alpha(V(t)) - \frac{5(H_{max}\bar{K}^2(t) + \varepsilon_1)}{\delta\mu\gamma}\right)\gamma\exp(\beta V(t))V(t) - \frac{\mu\delta}{4H_{max}}\varphi(V(t))\|v_x[t]\|_2^2$$
$$-\left(\alpha(V(t))\frac{\beta}{2} + \frac{\mu\delta\beta\gamma}{4H_{max}}\varphi(V(t)) - \frac{\varepsilon_2}{p_1^2(V(t))}\right)V(t)\exp(\beta V(t))\|v_x[t]\|_2^2$$
$$-\alpha(V(t))\beta\gamma V^2(t)\exp(\beta V(t)) - \left(\frac{\mu\delta\beta}{8H_{max}}\varphi(V(t)) - \frac{5L}{6\mu}\right)\exp(\beta V(t))\|v_x[t]\|_2^4$$
$$-\left(\delta + 1 - \frac{\mu\bar{K}(t)}{2g}\right)\left(\frac{\beta}{2}\|v_x[t]\|_2^2 + \gamma + \beta\gamma V(t)\right)\exp(\beta V(t))\int_0^L \kappa(h(t,x), v(t,x))v^2(t,x)dx$$
$$-\left(\delta + 1 - \frac{\mu\bar{K}(t)}{2g}\right)\int_0^L \kappa(h(t,x), v(t,x))v^2(t,x)dx$$

(3.49)

where $\bar{K}(t)$ is defined by (3.46), $V(t) = V(\xi(t), w(t), h[t], v[t])$ and $\varepsilon_1, \varepsilon_2 > 0$ are defined by (3.37), (3.38), respectively.

## 4. Proofs

**Proof of Lemma 1:** Notice that definition (3.6) implies that the equation $G(h) = \int_{h^*}^h \frac{|r - h^*|}{\sqrt{r}}dr$ holds for all $h > 0$.

Let $(\xi, w, h, v) \in S$ be given. Using the inequality

$$(h(x)v(x) + \mu h'(x))^2 \geq \frac{\delta\mu^2}{1+\delta}(h'(x))^2 - \delta h^2(x)v^2(x)$$

and definition (3.4), we obtain the following estimate:

$$W(h, v) \geq \frac{\delta\mu^2}{2(\delta+1)}\int_0^L h^{-1}(x)(h'(x))^2 dx - \frac{\delta}{2}\int_0^L h(x)v^2(x)dx + \frac{1}{2}g\|h - h^*\chi_{[0,L]}\|_2^2 \quad (4.1)$$

Using definitions (3.2), (3.3), (3.4) and estimate (4.1), we obtain:

$$\frac{\delta\mu^2}{2(\delta+1)}\int_0^L h^{-1}(x)(h'(x))^2 dx + \frac{1+\delta}{2}g\|h - h^*\chi_{[0,L]}\|_2^2 \leq V(\xi, w, h, v) \quad (4.2)$$

Let arbitrary $x, y \in [0, L]$ be given. Using the Cauchy-Schwarz inequality, the fact that $\int_0^L h(x)dx = m$ (recall definition (3.1)) and the fact that $G'(h) = \dfrac{|h - h^*|}{\sqrt{h}}$ for all $h > 0$, we get:

$$|G(h(x)) - G(h(y))| \leq \left|\int_y^x G'(h(s))h'(s)ds\right| \leq \int_{\min(x,y)}^{\max(x,y)} |G'(h(s))h'(s)|ds$$

$$\leq \left(\int_{\min(x,y)}^{\max(x,y)} h^{-1}(s)(h'(s))^2 ds\right)^{1/2} \left(\int_{\min(x,y)}^{\max(x,y)} h(s)(G'(h(s)))^2 ds\right)^{1/2}$$

$$\leq \left(\int_0^L h^{-1}(s)(h'(s))^2 ds\right)^{1/2} \left(\int_0^L h(s)(G'(h(s)))^2 ds\right)^{1/2} \quad (4.3)$$

$$= \left(\int_0^L h^{-1}(s)(h'(s))^2 ds\right)^{1/2} \|h - h^* \chi_{[0,L]}\|_2$$

$$|h(x) - h(y)| \leq \left|\int_y^x h'(s)ds\right| \leq \int_{\min(x,y)}^{\max(x,y)} |h'(s)|ds$$

$$\leq \left(\int_{\min(x,y)}^{\max(x,y)} h^{-1}(s)(h'(s))^2 ds\right)^{1/2} \left(\int_{\min(x,y)}^{\max(x,y)} h(s)ds\right)^{1/2} \quad (4.4)$$

$$\leq \left(\int_0^L h^{-1}(s)(h'(s))^2 ds\right)^{1/2} \left(\int_0^L h(s)ds\right)^{1/2} = m^{1/2}\left(\int_0^L h^{-1}(s)(h'(s))^2 ds\right)^{1/2}$$

Since $\max\left\{ab : \dfrac{\delta \mu^2}{2(\delta+1)}a^2 + \dfrac{(1+\delta)g}{2}b^2 \leq V\right\} = cV$ for all $V \geq 0$, we obtain from (4.2), (4.3), (4.4):

$$|G(h(x)) - G(h(y))| \leq cV(\xi, w, h, v), \text{ for all } x, y \in [0, L] \quad (4.5)$$

$$|h(x) - h(y)| \leq \frac{1}{\mu}\sqrt{\frac{2m(1+\delta)}{\delta}} V^{1/2}(\xi, w, h, v), \text{ for all } x, y \in [0, L] \quad (4.6)$$

Since $\int_0^L h(x)dx = m$ (recall definition (4.1)) and since $h^* = m/L$, it follows from continuity of $h$ and the mean value theorem that there exists $x^* \in [0, L]$ such that $h(x^*) = h^*$. Moreover, since $G(h^*) = 0$ (a consequence of definition (3.6)), we get from (4.5), (4.6) (with $y = x^*$):

$$-cV(\xi, w, h, v) \leq G(h(x)) \leq cV(\xi, w, h, v), \text{ for all } x \in [0, L] \quad (4.7)$$

$$h^* - \frac{1}{\mu}\sqrt{\frac{2m(1+\delta)}{\delta}} \sqrt{V(\xi, w, h, v)} \leq h(x) \leq h^* + \frac{1}{\mu}\sqrt{\frac{2m(1+\delta)}{\delta}} \sqrt{V(\xi, w, h, v)},$$
$$\text{for all } x \in [0, L] \quad (4.8)$$

Inequality (3.8) is a direct consequence of definitions (3.9) and estimates (4.7), (4.8). The proof is complete. ◁

**Proof of Lemma 2:** Using (3.3) we obtain for all $t > 0$:

$$\frac{d}{dt}E(h[t], v[t]) := \frac{1}{2}\int_0^L h_t(t,x)v^2(t,x)dx + \int_0^L h(t,x)v(t,x)v_t(t,x)dx \\ + g\int_0^L h_t(t,x)\left(h(t,x) - h^*\right)dx \quad (4.9)$$

Using (2.8) and (2.9) we obtain

$$v_t + vv_x + gh_x = \mu h^{-1}(hv_x)_x - h^{-1}k(h,v)v + f \quad (4.10)$$

Combining (2.8), (4.9) and (4.10) we get for all $t > 0$:

$$\frac{d}{dt}E(h[t], v[t]) = -\frac{1}{2}\int_0^L (h(t,x)v(t,x))_x v^2(t,x)dx - \int_0^L h(t,x)v^2(t,x)v_x(t,x)dx \\ - g\int_0^L h(t,x)v(x)h_x(t,x)dx + \mu\int_0^L v(t,x)(h(t,x)v_x(t,x))_x dx \\ + f(t)\int_0^L h(t,x)v(t,x)dx - \int_0^L \kappa(h(t,x), v(t,x))v^2(t,x)dx \\ - g\int_0^L \left(h(t,x) - h^*\right)(h(t,x)v(t,x))_x dx \quad (4.11)$$

Integrating by parts and using (2.10), we get for all $t > 0$

$$-g\int_0^L \left(h(t,x) - h^*\right)(h(t,x)v(t,x))_x dx = g\int_0^L h_x(t,x)h(t,x)v(t,x)dx \\ \mu\int_0^L v(t,x)(h(t,x)v_x(t,x))_x dx = -\mu\int_0^L h(t,x)v_x^2(t,x)dx \quad (4.12) \\ -\frac{1}{2}\int_0^L (h(t,x)v(t,x))_x v^2(t,x)dx = \int_0^L h(t,x)v^2(t,x)v_x(t,x)dx$$

As a consequence of (4.11) and (4.12) equation (3.40) holds for all $t > 0$.

Now we define for all $t \geq 0$ and $x \in [0, L]$ the following:

$$\varphi(t,x) := h(t,x)v(t,x) + \mu h_x(t,x) \quad (4.13)$$

Using definition (4.13) with (2.8) and (2.9) we get for all $t > 0$ and $x \in (0, L)$:

$$\varphi_t(t,x) = -\left(\varphi(t,x)v(t,x) + \frac{1}{2}gh^2(t,x)\right)_x - \kappa(h(t,x), v(t,x))v(t,x) + h(t,x)f(t) \quad (4.14)$$

Using (3.4) and definition (4.13), we obtain all $t > 0$:

$$\frac{d}{dt}W(h[t],v[t]) := -\frac{1}{2}\int_0^L h^{-2}(t,x)h_t(t,x)\varphi^2(t,x)dx + \int_0^L h^{-1}(t,x)\varphi(t,x)\varphi_t(t,x)dx$$
$$+ g\int_0^L \big(h(t,x)-h^*\big)h_t(t,x)dx \tag{4.15}$$

Combining (2.8), (4.14) and (4.15), we get for all $t > 0$:

$$\frac{d}{dt}W(h[t],v[t]) := \frac{1}{2}\int_0^L h^{-2}(t,x)\varphi^2(t,x)\big(h(t,x)v(t,x)\big)_x dx$$
$$-\int_0^L h^{-1}(t,x)\varphi(t,x)\left(\varphi(t,x)v(t,x)+\frac{1}{2}gh^2(t,x)\right)_x dx$$
$$-\int_0^L h^{-1}(t,x)\varphi(t,x)\kappa(h(t,x),v(t,x))v(t,x)dx \tag{4.16}$$
$$+ f(t)\int_0^L \varphi(t,x)dx - g\int_0^L \big(h(t,x)-h^*\big)\big(h(t,x)v(t,x)\big)_x dx$$

Using (2.10) and integration by parts, we get the following equation:

$$\frac{1}{2}\int_0^L h^{-2}(t,x)\varphi^2(t,x)\big(h(t,x)v(t,x)\big)_x dx =$$
$$\int_0^L h^{-2}(t,x)h_x(t,x)\varphi^2(t,x)v(t,x)dx - \int_0^L h^{-1}(t,x)\varphi(t,x)\varphi_x(t,x)v(t,x)dx \tag{4.17}$$

From (4.16) and (4.17) we conclude that

$$\frac{d}{dt}W(h[t],v[t]) = \int_0^L h^{-2}(t,x)h_x(t,x)\varphi^2(t,x)v(t,x)dx$$
$$-2\int_0^L h^{-1}(t,x)\varphi(t,x)\varphi_x(t,x)v(t,x)dx - \int_0^L h^{-1}(x)\varphi^2(t,x)v_x(t,x)dx$$
$$-\int_0^L h^{-1}(x)\varphi(t,x)\kappa(h(t,x),v(t,x))v(t,x)dx + f(t)\int_0^L \varphi(t,x)dx \tag{4.18}$$
$$+ g\int_0^L \big(h(t,x)v(t,x)-\varphi(t,x)\big)h_x(t,x)dx$$

for all $t > 0$. Equation (4.18) and definition (4.13) allow us to get for all $t > 0$:

$$\frac{d}{dt}W(h[t],v[t]) = -\int_0^L \left(h^{-1}(t,x)\big(h(t,x)v(t,x)+\mu h_x(t,x)\big)^2 v(t,x)\right)_x dx$$
$$-\int_0^L \kappa(h(t,x),v(t,x))v^2(t,x)dx - \mu\int_0^L h^{-1}(x)h_x(t,x)\kappa(h(t,x),v(t,x))v(t,x)dx \tag{4.19}$$
$$+ f(t)\int_0^L \big(h(t,x)v(t,x)+\mu h_x(t,x)\big)dx - g\mu\|h_x[t]\|_2^2$$

Equation (3.41) is a consequence of (2.10) and (4.19). ◁

**Proof of Lemma 3:** Let arbitrary $(\xi, w, h, v) \in X$ with $v \in H^1(0, L)$ and $V(\xi, w, h, v) < R$ be given. Due to the fact that $\int_0^L h(x)dx = m$ (recall (2.12)), the continuity of $h$ and the mean value theorem, there exists $x^* \in [0, L]$ such that $h(x^*) = h^* = m/L$. Using the Cauchy-Schwarz inequality, we get for all $x \in [0, L]$:

$$\left|h(x) - h^*\right| = \left|\int_{x^*}^x h'(s)ds\right| \leq \int_{\min(x^*,x)}^{\max(x,x^*)} |h'(s)|ds \leq \int_0^L |h'(s)|ds \leq \sqrt{L}\|h'\|_2 \quad (4.20)$$

Thus

$$\left(h(x) - h^*\right)^2 \leq L\|h'\|_2^2 \quad (4.21)$$

Using the inequality

$$\left(h(x)v(x) + \mu h'(x)\right)^2 \leq 2h^2(x)v^2(x) + 2\mu^2 \left(h'(x)\right)^2$$

and (3.8) we obtain

$$\int_0^L h^{-1}(x)\left(h(x)v(x) + \mu h'(x)\right)^2 dx \leq 2\int_0^L h(x)v^2(x)dx + \frac{2\mu^2}{p_1(V(\xi, w, h, v))}\|h'\|_2^2 \quad (4.22)$$

Since $v(0) = v(L) = 0$ by virtue of Wirtinger's inequality and (3.8) we have:

$$\|v\|_2^2 \leq \frac{L^2}{\pi^2}\|v'\|_2^2 \leq \frac{L^2}{\pi^2 p_1(V(\xi, w, h, v))}\int_0^L h(x)\left(v'(x)\right)^2 dx$$

Combining the above inequality and (3.8) we get

$$\int_0^L h(x)v^2(x)dx \leq p_2(V(\xi, w, h, v))\|v\|_2^2 \leq \frac{L^2}{\pi^2}\frac{p_2(V(\xi, w, h, v))}{p_1(V(\xi, w, h, v))}\int_0^L h(x)\left(v'(x)\right)^2 dx \quad (4.23)$$

Using (4.21), (4.22) and (4.23) we obtain

$$V(\xi, w, h, v) := \frac{\delta}{2}\int_0^L h(x)v^2(x)dx + \frac{\delta+1}{2}g\left\|h - h^*\chi_{[0,L]}\right\|_2^2$$
$$+\frac{1}{2}\int_0^L h^{-1}(x)\left(h(x)v(x) + \mu h'(x)\right)^2 dx + \frac{qk^2}{2}\xi^2 + \frac{q}{2}(w + k\xi)^2$$
$$\leq \frac{\delta+2}{2}\int_0^L h(x)v^2(x)dx + \left(\frac{\mu^2}{p_1(V(\xi, w, h, v))} + \frac{\delta+1}{2}gL^2\right)\|h'\|_2^2 + \frac{qk^2}{2}\xi^2 + \frac{q}{2}(w + k\xi)^2$$
$$\leq \Lambda(V(\xi, w, h, v))\left(\int_0^L h(x)\left(v'(x)\right)^2 dx + \|h'\|_2^2 + \xi^2 + (w + k\xi)^2\right)$$

(4.24)

where

$$\Lambda(s) := \frac{1}{2}\max\left(\frac{L^2(\delta+2)p_2(s)}{\pi^2 p_1(s)}, (\delta+1)gL^2 + \frac{2\mu^2}{p_1(s)}, qk^2, q\right), \text{ for } s \in [0, R) \quad (4.25)$$

Since $p_2 : \Re_+ \to \Re$ is an increasing function and $p_1 : \Re_+ \to \Re$ is a decreasing function, it follows from definition (4.25) that $\Lambda : [0, R) \to (0, +\infty)$ is a non-decreasing function.
The proof is complete. ◁

**Proof of Lemma 4:** Let arbitrary $(\xi, w, h, v) \in X$ with $V(\xi, w, h, v) < R$ be given. Using definitions (3.2), (3.3), (3.4), the inequalities

$$(h(x)v(x) + \mu h'(x))^2 \leq 2h^2(x)v^2(x) + 2\mu^2 (h'(x))^2,$$
$$(w + k\xi)^2 \leq 2w^2 + 2k^2\xi^2$$

and (3.12) we obtain

$$\begin{aligned}V(\xi, w, h, v) &\leq \frac{\delta + 2}{2} \int_0^L h(x)v^2(x)dx + \frac{\delta + 1}{2} g \left\| h - h^* \chi_{[0,L]} \right\|_2^2 \\ &+ \mu^2 \int_0^L h^{-1}(x)(h'(x))^2 dx + \frac{3qk^2}{2}\xi^2 + qw^2 \\ &\leq \frac{\delta + 2}{2} H_{\max} \|v\|_2^2 + \frac{\delta + 1}{2} g \left\| h - h^* \chi_{[0,L]} \right\|_2^2 \\ &+ \frac{\mu^2}{p_1(V(\xi, w, h, v))} \|h'\|_2^2 + \frac{3qk^2}{2}\xi^2 + qw^2\end{aligned} \quad (4.26)$$

The above inequality implies the left inequality (3.43) with

$$G_2(s) := \max\left(\frac{\delta + 2}{2} H_{\max}, \frac{\delta + 1}{2} g, \frac{\mu^2}{p_1(s)}, \frac{3qk^2}{2}, q\right), \text{ for } s \in [0, R)$$

Since $p_1 : \mathfrak{R}_+ \to \mathfrak{R}$ is a decreasing function, it follows from the above definition that $G_2 : [0, R) \to (0, +\infty)$ is a non-decreasing function.

Using definitions (3.2), (3.3), (3.4) and the inequalities

$$v(x)h'(x) \geq -\frac{\delta + 2}{4\mu} h(x)v^2(x) - \frac{\mu}{\delta + 2} h^{-1}(x)(h'(x))^2,$$
$$w\xi \geq -\frac{3k}{4}\xi^2 - \frac{1}{3k}w^2$$

we obtain

$$\begin{aligned}V(\xi, w, h, v) &\geq \frac{\delta}{4} \int_0^L h(x)v^2(x)dx + \frac{\delta + 1}{2} g \left\| h - h^* \chi_{[0,L]} \right\|_2^2 \\ &+ \frac{\delta\mu^2}{2\delta + 4} \int_0^L h^{-1}(x)(h'(x))^2 dx + \frac{qk^2}{4}\xi^2 + \frac{q}{6}w^2 \\ &\geq \frac{\delta}{4} p_1(V(\xi, w, h, v)) \|v\|_2^2 + \frac{\delta + 1}{2} g \left\| h - h^* \chi_{[0,L]} \right\|_2^2 \\ &+ \frac{\delta\mu^2}{H_{\max}(2\delta + 4)} \|h'\|_2^2 + \frac{qk^2}{4}\xi^2 + \frac{q}{6}w^2\end{aligned} \quad (4.27)$$

Inequality (4.27) implies the right inequality of (3.43) with

$$G_1(s) := \frac{2}{\min\left(\frac{\delta}{2} p_1(s), g(\delta + 1), \frac{\delta\mu^2}{H_{\max}(\delta + 2)}, \frac{qk^2}{2}, \frac{q}{3}\right)}, \text{ for } s \in [0, R)$$

Since $p_1 : \mathfrak{R}_+ \to \mathfrak{R}$ is a decreasing function, it follows from the above definition that $G_1 : [0, R) \to (0, +\infty)$ is a non-decreasing function.
The proof is complete. ◁

**Proof of Proposition 1:** Suppose first that $L = 1$. Let $\varphi \in H^1(0,1)$ with $\varphi(0) = \varphi(1) = 0$ be a given function. It holds that $\varphi = \sum_{n=1}^{\infty} a_n \phi_n$, where $\phi_n(x) = \sqrt{2} \sin(n\pi x)$ and $a_n = \int_0^1 \varphi(x) \phi_n(x) dx$. Thus, we have (using the Cauchy-Schwarz inequality):

$$\|\varphi\|_\infty \leq \sqrt{2} \sum_{n=1}^{\infty} |a_n| = \sqrt{2} \sum_{n=1}^{\infty} \frac{1}{n\pi} n\pi |a_n| \leq \sqrt{2} \left( \sum_{n=1}^{\infty} \frac{1}{n^2 \pi^2} \right)^{1/2} \left( \sum_{n=1}^{\infty} n^2 \pi^2 a_n^2 \right)^{1/2} \quad (4.28)$$

Since $\varphi' = \sum_{n=1}^{\infty} n\pi a_n \psi_n$, where $\psi_n(x) = \sqrt{2} \cos(n\pi x)$ we obtain from Parseval's identity that

$$\|\varphi'\|_2^2 = \sum_{n=1}^{\infty} n^2 \pi^2 a_n^2 \quad (4.29)$$

Inequality (3.19) for $L = 1$ is a consequence of (4.28), (4.29) and the fact that $\sum_{n=1}^{\infty} \frac{1}{n^2 \pi^2} = \frac{1}{6}$. Since $\varphi'' = -\sum_{n=1}^{\infty} n^2 \pi^2 a_n \phi_n$ when $\varphi \in H^2(0,1)$, we have from Parseval's identity the equation $\|\varphi''\|_2^2 = \sum_{n=1}^{\infty} n^4 \pi^4 a_n^2$. Thus, (3.20) for $L = 1$ follows from (4.29) and the previous equation.

When $L \neq 1$, we consider the function $\tilde{\varphi}(\xi) = \varphi(L\xi)$ defined for $\xi \in [0,1]$.
The proof is complete. ◁

**Proof of Lemma 5:** Let $\sigma, k, q, \delta > 0$ be given constants and let $r \in [0, R)$ be a constant, where $R > 0$ is defined by (3.11). Consider a classical solution of system (2.7)-(2.10), (2.12) and (3.27) with $v \in C^0(\mathfrak{R}_+; H^1(0,L)) \cap C^1((0,+\infty); H^1(0,L))$ at a time $t > 0$ for which $V(\xi(t), w(t), h[t], v[t]) < R$. Using Lemma 2, (3.40), (3.41) and definition (3.2) we obtain:

$$\frac{d}{dt} V(\xi(t), w(t), h[t], v[t])$$

$$= -\mu g \|h_x[t]\|_2^2 - \mu \delta \int_0^L h(t,x) v_x^2(t,x) dx + qk(w(t) + k\xi(t))^2$$

$$-qk^3 \xi^2(t) + f(t) \left( (\delta + 1) \int_0^L h(t,x) v(t,x) dx + \mu(h(t,L) - h(t,0)) - q(w(t) + k\xi(t)) \right)$$

$$-(\delta + 1) \int_0^L \kappa(h(t,x), v(t,x)) v^2(t,x) dx - \mu \int_0^L h^{-1}(t,x) h_x(t,x) \kappa(h(t,x), v(t,x)) v(t,x) dx$$

(4.30)

Using (3.27) and (4.30), we get:

$$\frac{d}{dt}V(\xi(t),w(t),h[t],v[t])$$

$$= -\mu g \|h_x[t]\|_2^2 - \mu\delta\int_0^L h(t,x)v_x^2(t,x)dx + qk(w(t)+k\xi(t))^2$$

$$-qk^3\xi^2(t) - \sigma\left((\delta+1)\int_0^L h(t,x)v(t,x)dx + \mu(h(t,L)-h(t,0)) - q(w(t)+k\xi(t))\right)^2$$

$$-(\delta+1)\int_0^L \kappa(h(t,x),v(t,x))v^2(t,x)dx - \mu\int_0^L h^{-1}(t,x)h_x(t,x)\kappa(h(t,x),v(t,x))v(t,x)dx$$

(4.31)

Using the inequality

$$2\left((\delta+1)\int_0^L h(t,x)v(t,x)dx + \mu(h(t,L)-h(t,0))\right)(w(t)+k\xi(t))$$

$$\leq \frac{1+z}{q}\left((\delta+1)\int_0^L h(t,x)v(t,x)dx + \mu(h(t,L)-h(t,0))\right)^2 + \frac{q}{1+z}(w(t)+k\xi(t))^2$$

with $z = \dfrac{g\mu\delta\pi^2 p_1(r)}{2\sigma L\left(mgLH_{\max}(\delta+1)^2 + 2\mu^2\delta\pi^2 p_1(r)\right)}$, we obtain from (4.31) the following estimate:

$$\frac{d}{dt}V(\xi(t),w(t),h[t],v[t]) \leq -\mu g \|h_x[t]\|_2^2$$

$$-\mu\delta\int_0^L h(t,x)v_x^2(t,x)dx - qk^3\xi^2(t) - q\left(\frac{\sigma qz}{1+z}-k\right)(w(t)+k\xi(t))^2$$

$$+\sigma z\left((\delta+1)\int_0^L h(t,x)v(t,x)dx + \mu(h(t,L)-h(t,0))\right)^2 \quad (4.32)$$

$$-(\delta+1)\int_0^L \kappa(h(t,x),v(t,x))v^2(t,x)dx$$

$$-\mu\int_0^L h^{-1}(t,x)h_x(t,x)\kappa(h(t,x),v(t,x))v(t,x)dx$$

Using the inequality

$$2\mu(\delta+1)\left(\int_0^L h(t,x)v(t,x)dx\right)(h(t,L)-h(t,0))$$

$$\leq \frac{2\mu^2\delta\pi^2 p_1(r)}{mgLH_{\max}}\left(\int_0^L h(t,x)v(t,x)dx\right)^2 + \frac{mgLH_{\max}(\delta+1)^2}{2\delta\pi^2 p_1(r)}\left(\int_0^L h_x(t,x)dx\right)^2$$

and the fact that $z = \dfrac{g\mu\delta\pi^2 p_1(r)}{2\sigma L\left(mgLH_{\max}(\delta+1)^2 + 2\mu^2\delta\pi^2 p_1(r)\right)}$ we obtain from (4.32) the following estimate:

$$\frac{d}{dt}V(\xi(t),w(t),h[t],v[t]) \leq -\mu g \|h_x[t]\|_2^2 - \mu\delta\int_0^L h(t,x)v_x^2(t,x)dx$$

$$-qk^3\xi^2(t) + \frac{\mu\delta\pi^2 p_1(r)}{2mL^2 H_{\max}}\left(\int_0^L h(t,x)v(t,x)dx\right)^2$$

$$-q(q\theta(r)-k)(w(t)+k\xi(t))^2 + \frac{g\mu}{4L}\left(\int_0^L h_x(t,x)dx\right)^2 \quad (4.33)$$

$$-(\delta+1)\int_0^L \kappa(h(t,x),v(t,x))v^2(t,x)dx$$

$$-\mu\int_0^L h^{-1}(t,x)h_x(t,x)\kappa(h(t,x),v(t,x))v(t,x)dx$$

where $\theta(r)$ is defined by (3.26). The Cauchy-Schwarz inequality in conjunction with (2.12) implies the following inequalities:

$$\left(\int_0^L h_x(t,x)dx\right)^2 \leq L\|h_x[t]\|_2^2$$

$$\left(\int_0^L h(t,x)v(t,x)dx\right)^2 \leq m\int_0^L h(t,x)v^2(t,x)dx$$

Thus, we obtain from (4.33) the following estimate:

$$\frac{d}{dt}V(\xi(t),w(t),h[t],v[t]) \leq -\frac{3\mu g}{4}\|h_x[t]\|_2^2$$

$$-\mu\delta\int_0^L h(t,x)v_x^2(t,x)dx - (\delta+1)\int_0^L \kappa(h(t,x),v(t,x))v^2(t,x)dx$$

$$-qk^3\xi^2(t) + \frac{\mu\delta\pi^2 p_1(r)}{2L^2 H_{\max}}\int_0^L h(t,x)v^2(t,x)dx - q(q\theta(r)-k)(w(t)+k\xi(t))^2 \quad (4.34)$$

$$-\mu\int_0^L h^{-1}(t,x)h_x(t,x)\kappa(h(t,x),v(t,x))v(t,x)dx$$

Using the inequality

$$h^{-1}(t,x)h_x(t,x)\kappa(h(t,x),v(t,x))v(t,x) \geq -\frac{g}{2}h_x^2(t,x) - \frac{\kappa^2(h(t,x),v(t,x))}{2gh^2(t,x)}v^2(t,x)$$

we obtain from (4.34) the following estimate:

$$\frac{d}{dt}V(\xi(t),w(t),h[t],v[t]) \leq -\frac{\mu g}{4}\|h_x[t]\|_2^2 - \mu\delta\int_0^L h(t,x)v_x^2(t,x)dx$$
$$-(\delta+1)\int_0^L\left(1-\frac{\mu\kappa(h(t,x),v(t,x))}{2g(\delta+1)h^2(t,x)}\right)\kappa(h(t,x),v(t,x))v^2(t,x)dx \qquad (4.35)$$
$$-qk^3\xi^2(t) + \frac{\mu\delta\pi^2 p_1(r)}{2L^2 H_{\max}}\int_0^L h(t,x)v^2(t,x)dx - q(q\theta(r)-k)(w(t)+k\xi(t))^2$$

Definition (3.46) in conjunction with (3.12) and (4.35) implies the following estimate:

$$\frac{d}{dt}V(\xi(t),w(t),h[t],v[t]) \leq -\frac{\mu g}{4}\|h_x[t]\|_2^2 - \mu\delta\int_0^L h(t,x)v_x^2(t,x)dx$$
$$-\left(\delta+1-\frac{\mu\overline{K}(t)}{2g}\right)\int_0^L \kappa(h(t,x),v(t,x))v^2(t,x)dx \qquad (4.36)$$
$$-qk^3\xi^2(t) + \frac{\mu\delta\pi^2 p_1(r)p_2(V(t))}{2L^2 H_{\max}}\|v[t]\|_2^2 - q(q\theta(r)-k)(w(t)+k\xi(t))^2$$

Since $v(t,0)=v(t,L)=0$ (recall (2.10)), by virtue of Wirtinger's inequality and (3.12), we get:

$$\|v[t]\|_2^2 \leq \frac{L^2}{\pi^2}\|v_x[t]\|_2^2 \leq \frac{L^2}{\pi^2 p_1(V(t))}\int_0^L h(t,x)v_x^2(t,x)dx \qquad (4.37)$$

Combining (4.36) and (4.37), we obtain (3.44).

Since $v \in C^0(\Re_+;H^1(0,L)) \cap C^1((0,+\infty);H^1(0,L))$, by working with a similar way with the proof of Lemma 3.3 in [20], we establish the following equation:

$$\frac{d}{dt}\left(\frac{1}{2}\|v_x[t]\|_2^2\right) = -\int_0^L v_{xx}(t,x)v_t(t,x)dx \qquad (4.38)$$

Using (2.8), (2.9) and (4.38) we get:

$$\frac{d}{dt}\left(\frac{1}{2}\|v_x[t]\|_2^2\right) = \int_0^L v_{xx}(t,x)v(t,x)v_x(t,x)dx + g\int_0^L v_{xx}(t,x)h_x(t,x)dx$$
$$-\mu\|v_{xx}[t]\|_2^2 - f(t)\int_0^L v_{xx}(t,x)dx - \mu\int_0^L v_{xx}(t,x)h^{-1}(t,x)h_x(t,x)v_x(t,x)dx \qquad (4.39)$$
$$+\int_0^L v_{xx}(t,x)h^{-1}(t,x)\kappa(h(t,x),v(y,x))v(t,x)dx$$

Using (4.39) and the following Young inequalities

$$v_{xx}(t,x)v(t,x)v_x(t,x) \leq \frac{\mu}{10}v_{xx}^2(t,x) + \frac{5}{2\mu}v^2(t,x)v_x^2(t,x)$$

$$gv_{xx}(t,x)h_x(t,x) \leq \frac{\mu}{10}v_{xx}^2(t,x) + \frac{5g^2}{2\mu}h_x^2(t,x)$$

$$-\mu v_{xx}(t,x)h^{-1}(t,x)h_x(t,x)v_x(t,x) \leq \frac{\mu}{10}v_{xx}^2(t,x) + \frac{5\mu}{2}h^{-2}(t,x)h_x^2(t,x)v_x^2(t,x)$$

$$-f(t)v_{xx}(t,x) \leq \frac{\mu}{10}v_{xx}^2(t,x) + \frac{5}{2\mu}f^2(t)$$

$$v_{xx}(t,x)h^{-1}(t,x)\kappa(h(t,x),v(t,x))v(t,x)$$
$$\leq \frac{\mu}{10}v_{xx}^2(t,x) + \frac{5}{2\mu}h^{-2}(t,x)\kappa^2(h(t,x),v(t,x))v^2(t,x)$$

we get:

$$\frac{d}{dt}\left(\frac{1}{2}\|v_x[t]\|_2^2\right) \leq \frac{5}{2\mu}\int_0^L v^2(t,x)v_x^2(t,x)dx + \frac{5g^2}{2\mu}\|h_x[t]\|_2^2 - \frac{\mu}{2}\|v_{xx}[t]\|_2^2 + \frac{5L}{2\mu}f^2(t)$$
$$+ \frac{5\mu}{2}\int_0^L h^{-2}(t,x)h_x^2(t,x)v_x^2(t,x)dx + \frac{5}{2\mu}\int_0^L h^{-2}(t,x)\kappa^2(h(t,x),v(t,x))v^2(t,x)dx \quad (4.40)$$

Applying Proposition 1 to the function $\varphi(\xi) = v(t,L\xi)$, we obtain:

$$\|v[t]\|_\infty \leq \sqrt{\frac{L}{3}}\|v_x[t]\|_2 \text{ and } \frac{\pi}{L}\|v_x[t]\|_2 \leq \|v_{xx}[t]\|_2 \quad (4.41)$$

Thus, we get from (4.41):

$$\int_0^L v^2(t,x)v_x^2(t,x)dx \leq \|v[t]\|_\infty^2 \|v_x[t]\|_2^2 \leq \frac{L}{3}\|v_x[t]\|_2^4 \quad (4.42)$$

Combining (4.40) and (4.42) we obtain:

$$\frac{d}{dt}\left(\frac{1}{2}\|v_x[t]\|_2^2\right) \leq \frac{5L}{6\mu}\|v_x[t]\|_2^4 + \frac{5g^2}{2\mu}\|h_x[t]\|_2^2 - \frac{\mu}{2}\|v_{xx}[t]\|_2^2 + \frac{5L}{2\mu}f^2(t)$$
$$+ \frac{5\mu}{2}\int_0^L h^{-2}(t,x)h_x^2(t,x)v_x^2(t,x)dx + \frac{5}{2\mu}\int_0^L h^{-2}(t,x)\kappa^2(h(t,x),v(t,x))v^2(t,x)dx \quad (4.43)$$

Using (3.12), (3.46) and (4.43) we get:

$$\frac{d}{dt}\left(\frac{1}{2}\|v_x[t]\|_2^2\right) \leq \frac{5L}{6\mu}\|v_x[t]\|_2^4 + \frac{5g^2}{2\mu}\|h_x[t]\|_2^2 - \frac{\mu}{2}\|v_{xx}[t]\|_2^2 + \frac{5L}{2\mu}f^2(t)$$
$$+ \frac{5\mu}{2p_1(V(t))}\int_0^L h^{-1}(t,x)h_x^2(t,x)v_x^2(t,x)dx + \frac{5\bar{K}^2(t)}{2\mu}\int_0^L h^2(t,x)v^2(t,x)dx \quad (4.44)$$

Since $v(t, L) - v(t, 0) = \int_0^L v_x(t, x)dx = 0$ and since $v[t] \in C^1([0, L])$, there exists $x^* \in [0, L]$ such that $v_x(t, x^*) = 0$. Using the inequality $2|v_x(t, s)v_{xx}(t, s)| \le \varepsilon v_x^2(t, s) + \varepsilon^{-1} v_{xx}^2(t, s)$ that holds for all $\varepsilon > 0$, we get for all $x \in [0, L]$:

$$v_x^2(t, x) = 2\int_{x^*}^{x} v_x(t, s)v_{xx}(t, s)ds \le 2 \int_{\min(x,x^*)}^{\max(x,x^*)} |v_x(t, s)v_{xx}(t, s)|ds$$

$$\le \varepsilon \int_{\min(x,x^*)}^{\max(x,x^*)} v_x^2(t, s)ds + \varepsilon^{-1} \int_{\min(x,x^*)}^{\max(x,x^*)} v_{xx}^2(t, s)ds \le \varepsilon \|v_x[t]\|_2^2 + \varepsilon^{-1} \|v_{xx}[t]\|_2^2$$

Thus, we obtain for all $\varepsilon > 0$:

$$\|v_x[t]\|_\infty^2 \le \varepsilon \|v_x[t]\|_2^2 + \varepsilon^{-1} \|v_{xx}[t]\|_2^2, \text{ for all } \varepsilon > 0 \quad (4.45)$$

Combining (4.44) and (4.45), we get for all $\varepsilon > 0$:

$$\frac{d}{dt}\left(\frac{1}{2}\|v_x[t]\|_2^2\right) \le \frac{5L}{6\mu}\|v_x[t]\|_2^4 + \frac{5g^2}{2\mu}\|h_x[t]\|_2^2 - \frac{\mu}{2}\|v_{xx}[t]\|_2^2 + \frac{5L}{2\mu}f^2(t)$$

$$+ \frac{5\mu\varepsilon\|v_x[t]\|_2^2}{2p_1(V(t))}\int_0^L h^{-1}(t, x)h_x^2(t, x)dx + \frac{5\mu\|v_{xx}[t]\|_2^2}{2\varepsilon p_1(V(t))}\int_0^L h^{-1}(t, x)h_x^2(t, x)dx \quad (4.46)$$

$$+ \frac{5\bar{K}^2(t)}{2\mu}\int_0^L h^2(t, x)v^2(t, x)dx$$

Using (3.12) and assuming $\|h_x[t]\|_2^2 > 0$, we select $\varepsilon = \frac{10}{p_1(V(t))}\int_0^L h^{-1}(t, x)h_x^2(t, x)dx$ and we get:

$$\frac{d}{dt}\left(\frac{1}{2}\|v_x[t]\|_2^2\right) \le \frac{5L}{6\mu}\|v_x[t]\|_2^4 + \frac{5g^2 H_{\max}}{2\mu}\int_0^L h^{-1}(t, x)h_x^2(t, x)dx - \frac{\mu}{4}\|v_{xx}[t]\|_2^2$$

$$+ \frac{5L}{2\mu}f^2(t) + \frac{25\mu\|v_x[t]\|_2^2}{p_1^2(V(t))}\left(\int_0^L h^{-1}(t, x)h_x^2(t, x)dx\right)^2 + \frac{5\bar{K}^2(t)}{2\mu}\int_0^L h^2(t, x)v^2(t, x)dx \quad (4.47)$$

Inequality (4.47) also holds when $\|h_x[t]\|_2^2 = 0$. Thus, we conclude that (4.47) holds in every case.

Using inequality (4.2), which implies $\int_0^L h^{-1}(t, x)h_x^2(t, x)dx \le \frac{2(\delta+1)}{\delta\mu^2}V(t)$, we get from (4.47):

$$\frac{d}{dt}\left(\frac{1}{2}\|v_x[t]\|_2^2\right) \le \frac{5L}{6\mu}\|v_x[t]\|_2^4 + \frac{5(\delta+1)g^2 H_{\max}}{\delta\mu^3}V(t) - \frac{\mu}{4}\|v_{xx}[t]\|_2^2$$

$$+ \frac{5L}{2\mu}f^2(t) + \frac{100(\delta+1)^2\|v_x[t]\|_2^2}{\delta^2\mu^3 p_1^2(V(t))}V^2(t) + \frac{5\bar{K}^2(t)}{2\mu}\int_0^L h^2(t, x)v^2(t, x)dx \quad (4.48)$$

Using (3.2), (3.3) (which both imply that $\int_0^L h(t, x)v^2(t, x)dx \le \frac{2}{\delta}V(t)$), (3.12), the fact that $V(t) < R$ and (4.48) we get:

$$\frac{d}{dt}\left(\frac{1}{2}\|v_x[t]\|_2^2\right) \leq \frac{5L}{6\mu}\|v_x[t]\|_2^4 - \frac{\mu}{4}\|v_{xx}[t]\|_2^2 + \frac{5L}{2\mu}f^2(t)$$
$$+ \frac{100(\delta+1)^2 R}{\delta^2 \mu^3 p_1^2(V(t))}\|v_x[t]\|_2^2 V(t) + \frac{5H_{\max}}{\mu\delta}\left(\bar{K}^2(t) + \frac{(\delta+1)g^2}{\mu^2}\right)V(t) \quad (4.49)$$

Equation (3.27) and the triangle inequality give:

$$|f(t)| \leq \sigma\left((\delta+1)\int_0^L h(t,x)|v(t,x)|dx + \mu\|h_x[t]\|_1 + q|w(t) + k\xi(t)|\right) \quad (4.50)$$

Using the Cauchy-Schwarz inequality and (2.12), we get:

$$\|h_x[t]\|_1 \leq \sqrt{m}\left(\int_0^L h^{-1}(t,x)h_x^2(t,x)dx\right)^{1/2}$$
$$\int_0^L h(t,x)|v(t,x)|dx \leq \sqrt{m}\left(\int_0^L h(t,x)v^2(t,x)dx\right)^{1/2} \quad (4.51)$$

Using (3.2), (3.3) (which both imply that $\int_0^L h(t,x)v^2(t,x)dx \leq \frac{2}{\delta}V(t)$ and $(w(t)+k\xi(t))^2 \leq \frac{2}{q}V(t)$) and (4.2) (which implies that $\int_0^L h^{-1}(t,x)h_x^2(t,x)dx \leq \frac{2(\delta+1)}{\delta\mu^2}V(t)$), we get from (4.50), (4.51):

$$f^2(t)$$
$$\leq 3\sigma^2\left((\delta+1)^2 m\int_0^L h(t,x)v^2(t,x)dx + \mu^2 m\int_0^L h^{-1}(t,x)h_x^2(t,x)dx + q^2(w(t)+k\xi(t))^2\right) \quad (4.52)$$
$$\leq 6\sigma^2\left(\frac{(\delta+1)(\delta+2)m}{\delta} + q\right)V(t)$$

Combining (4.49) and (4.52) and using (4.41) we obtain (3.45).

The proof is complete. ◁

**Proof of Lemma 6:** Let $\beta, \gamma, \sigma, k, q, \delta > 0$ be given constants. Suppose that (3.25) holds. Consider a classical solution of the PDE-ODE system (2.7)-(2.10), (2.12) and (3.27) with $v \in C^0(\Re_+; H^1(0,L)) \cap C^1((0,+\infty); H^1(0,L))$ at a time $t > 0$ for which $V(\xi(t), w(t), h[t], v[t]) < R$ and $\varphi(V(\xi(t), w(t), h[t], v[t])) > 0$. Then we get from (3.44), (3.47) and (3.12):

$$\frac{d}{dt}V(\xi(t),w(t),h[t],v[t]) \le -\frac{\mu g}{4}\|h_x[t]\|_2^2 - q(q\theta(r)-k)(w(t)+k\xi(t))^2$$
$$-\frac{\mu\delta}{4H_{max}}\left(2H_{max}-p_1(r)\frac{p_2(V(t))}{p_1(V(t))}\right)\int_0^L h(t,x)v_x^2(t,x)dx - \frac{\mu\delta}{4H_{max}}\varphi(V(t))\|v_x[t]\|_2^2 \quad (4.53)$$
$$-\left(\delta+1-\frac{\mu\bar{K}(t)}{2g}\right)\int_0^L \kappa(h(t,x),v(t,x))v^2(t,x)dx - qk^3\xi^2(t)$$

Using definition (3.48), we obtain from (4.53):

$$\frac{d}{dt}V(\xi(t),w(t),h[t],v[t])$$
$$\le -\alpha(V(t))\Lambda(V(t))\left(\|h_x[t]\|_2^2 + (w(t)+k\xi(t))^2 + \xi^2(t) + \int_0^L h(t,x)v_x^2(t,x)dx\right) \quad (4.54)$$
$$-\frac{\mu\delta}{4H_{max}}\varphi(V(t))\|v_x[t]\|_2^2 - \left(\delta+1-\frac{\mu\bar{K}(t)}{2g}\right)\int_0^L \kappa(h(t,x),v(t,x))v^2(t,x)dx$$

Using (3.42) and (4.54) we obtain:

$$\frac{d}{dt}V(\xi(t),w(t),h[t],v[t]) \le -\alpha(V(t))V(t) - \frac{\mu\delta}{4H_{max}}\varphi(V(t))\|v_x[t]\|_2^2$$
$$-\left(\delta+1-\frac{\mu\bar{K}(t)}{2g}\right)\int_0^L \kappa(h(t,x),v(t,x))v^2(t,x)dx \quad (4.55)$$

Inequality (3.49) is a direct consequence of (4.55), (3.45) and definition (3.5). The proof is complete. ◁

**Proof of Theorem 1:** Let $r\in[0,R)$ be given and let constants $\sigma,q,k,\delta>0$, for which (3.24), (3.25) hold. Consider a classical solution of the PDE-ODE system (2.7)-(2.10), (2.12) with (3.27) with $v\in C^0(\Re_+;H^1(0,L))\cap C^1((0,+\infty);H^1(0,L))$ that satisfies $V(\xi(0),w(0),h[0],v[0])\le r$. Let $\bar{r}\in(r,R)$ be a constant that satisfies:

$$\max\left(\frac{p_2(\bar{r})}{p_1(\bar{r})} - \frac{2H_{max}}{p_1(r)}, K(p_1(\bar{r})) - \frac{2g(\delta+1)}{\mu}\right) < 0 \quad (4.56)$$

Such a constant $\bar{r}\in(r,R)$ exists because (3.24) holds and because all functions involved in (4.56) are continuous.

We next notice that the following implication holds:

"If $t>0$ and $V(\xi(t),w(t),h[t],v[t])\le\bar{r}$ then $\frac{d}{dt}V(\xi(t),w(t),h[t],v[t])\le 0$" \quad (4.57)

Indeed, implication (4.57) is a consequence of (3.44), (3.12), (3.23), (4.56) and definition (3.46).

We next prove by contradiction that $V(\xi(t), w(t), h[t], v[t]) \leq \bar{r}$ for all $t \geq 0$. Assume the contrary, i.e. that there exists $t \geq 0$ such that $V(\xi(t), w(t), h[t], v[t]) > \bar{r}$. Therefore, the set $A := \{t \geq 0 : V(\xi(t), w(t), h[t], v[t]) > \bar{r}\}$ is non-empty and bounded from below. Thus, we can define $t^* = \inf(A)$. By definition, and since $V(\xi(0), w(0), h[0], v[0]) \leq r < \bar{r}$ it holds that $t^* > 0$ and $V(\xi(t), w(t), h[t], v[t]) \leq \bar{r}$ for all $t \in [0, t^*)$. By continuity of the mapping $t \to V(\xi(t), w(t), h[t], v[t])$ (recall that $h \in C^1(\Re_+ \times [0, L]; (0, +\infty))$ and $v \in C^0(\Re_+ \times [0, L])$) we obtain that $V(\xi(t^*), w(t^*), h[t^*], v[t^*]) = \bar{r}$. Moreover, since $V(\xi(t), w(t), h[t], v[t]) \leq \bar{r}$ for all $t \in [0, t^*)$, it follows from implication (4.57) that $\frac{d}{dt} V(\xi(t), w(t), h[t], v[t]) \leq 0$ for all $t \in [0, t^*)$. By continuity of the mapping $t \to V(\xi(t), w(t), h[t], v[t])$, we obtain that

$$V(\xi(t^*), w(t^*), h[t^*], v[t^*]) \leq V(\xi(0), w(0), h[0], v[0]) \leq r < \bar{r}$$

which contradicts the fact that $V(\xi(t^*), w(t^*), h[t^*], v[t^*]) = \bar{r}$.

Since $V(\xi(t), w(t), h[t], v[t]) \leq \bar{r}$ for all $t \geq 0$, we conclude from implication (4.57) that $\frac{d}{dt} V(\xi(t), w(t), h[t], v[t]) \leq 0$ for all $t > 0$. By continuity of the mapping $t \to V(\xi(t), w(t), h[t], v[t])$, we obtain that

$$V(\xi(t), w(t), h[t], v[t]) \leq V(\xi(0), w(0), h[0], v[0]) \leq r < R, \text{ for all } t \geq 0 \quad (4.58)$$

Hence, $(\xi(t), w(t), h[t], v[t]) \in X_V(r)$ for all $t \geq 0$ (recall (3.16)). Using (3.44), (3.12), (3.23), definition (3.47) and (4.58), we get for all $t > 0$:

$$\frac{d}{dt} V(\xi(t), w(t), h[t], v[t]) \leq -\frac{\mu g}{4} \|h_x[t]\|_2^2 - q(q\theta(r) - k)(w(t) + k\xi(t))^2$$
$$- \frac{\mu \delta \varphi(r)}{2 H_{max} p_1(r)} \int_0^L h(t,x) v_x^2(t,x) dx - qk^3 \xi^2(t) \quad (4.59)$$

Thus, we obtain from (4.59) the following estimate for $t > 0$

$$\frac{d}{dt} V(\xi(t), w(t), h[t], v[t])$$
$$\leq -\omega \left( \int_0^L h_x^2(t,x) dx + \int_0^L h(t,x) v_x^2(t,x) dx + \xi^2(t) + (w(t) + k\xi(t))^2 \right) \quad (4.60)$$

where

$$\omega := \min\left( \frac{\mu g}{4}, \frac{\mu \delta \varphi(r)}{2 H_{max} p_1(r)}, qk^3, q(q\theta(r) - k) \right) \quad (4.61)$$

It follows from Lemma 3 and (4.60) that the following estimate holds for all $t > 0$:

$$\frac{d}{dt}V(\xi(t),w(t),h[t],v[t]) \leq -\frac{\omega V(\xi(t),w(t),h[t],v[t])}{\Lambda\big(V(\xi(t),w(t),h[t],v[t])\big)} \tag{4.62}$$

where $\Lambda$ is the non-decreasing function involved in (3.42). Since $\Lambda$ is non-decreasing, we obtain from (4.62) and (4.58) the following estimate for all $t>0$:

$$\frac{d}{dt}V(\xi(t),w(t),h[t],v[t]) \leq -\frac{\omega}{\Lambda(r)}V(\xi(t),w(t),h[t],v[t]) \tag{4.63}$$

Pick constants $\beta,\gamma>0$ such that

$$\begin{aligned}\beta\gamma &\geq \frac{4H_{max}\varepsilon_2}{\mu\delta p_1^2(r)\varphi(r)} \\ \beta &\geq \frac{20LH_{max}}{3\mu^2\delta\varphi(r)} \\ \gamma &> \frac{5\big(H_{max}K^2(p_1(r))+\varepsilon_1\big)}{\delta\mu\alpha(r)}\end{aligned} \tag{4.64}$$

We obtain from (3.46), (3.12), (3.23), (3.47), (3.48), (3.49) and (4.58) for all $t>0$:

$$\begin{aligned}\frac{d}{dt}U(\xi(t),w(t),h[t],v[t]) &\leq -\alpha(r)V(t) \\ &-\left(\frac{\delta\gamma\varphi(r)}{H_{max}}+\frac{\pi^2}{L^2}\right)\frac{\mu}{4}\exp(\beta V(t))\|v_x[t]\|_2^2 \\ &-\left(\alpha(r)-\frac{5\big(H_{max}K^2(p_1(r))+\varepsilon_1\big)}{\delta\mu\gamma}\right)\gamma\exp(\beta V(t))V(t)\end{aligned} \tag{4.65}$$

Thus, we obtain from (4.65) and (3.5) the following estimate for $t>0$:

$$\frac{d}{dt}U(\xi(t),w(t),h[t],v[t]) \leq -\bar{\omega}U(\xi(t),w(t),h[t],v[t]) \tag{4.66}$$

where

$$\bar{\omega}:=\min\left(\left(\frac{\delta\gamma\varphi(r)}{H_{max}}+\frac{\pi^2}{L^2}\right)\frac{\mu}{2},\alpha(r)-\frac{5\big(H_{max}K^2(p_1(r))+\varepsilon_1\big)}{\delta\mu\gamma}\right) \tag{4.67}$$

By continuity of the mappings $t\to V(\xi(t),w(t),h[t],v[t])$, $t\to U(\xi(t),w(t),h[t],v[t])$, (recall that $v\in C^0\big(\Re_+;H^1(0,L)\big)$, $h\in C^1\big(\Re_+\times[0,L];(0,+\infty)\big)$ and $v\in C^0\big(\Re_+\times[0,L]\big)$) the differential inequalities (4.63), (4.66) imply the following estimates for all $t\geq 0$:

$$V(\xi(t),w(t),h[t],v[t]) \leq \exp\left(-\frac{\omega t}{\Lambda(r)}\right)V(\xi(0),w(0),h[0],v[0]) \tag{4.68}$$

$$U(\xi(t),w(t),h[t],v[t]) \leq \exp(-\bar{\omega}t)U(\xi(0),w(0),h[0],v[0]) \tag{4.69}$$

Estimates (4.58), (4.68), (4.69), definition (3.5) and Lemma 4 give the following estimate for all $t \geq 0$:

$$\left\|(\xi(t), w(t), h[t] - h^* \chi_{[0,L]}, v[t])\right\|_X^2 \leq \Omega(t) \exp\left(-\frac{\omega t}{\Lambda(r)}\right) \left\|(\xi(0), w(0), h[0] - h^* \chi_{[0,L]}, v[0])\right\|_X^2 \tag{4.70}$$

$$\left\|v_x[t]\right\|_2^2 \leq \bar{\Omega} \exp(-\bar{\omega} t) \left( \left\|(\xi(0), w(0), h[0] - h^* \chi_{[0,L]}, v[0])\right\|_X^2 + \left\|v_x[0]\right\|_2^2 \right) \tag{4.71}$$

with

$$\Omega(t) := G_1\left(V\left(\xi(t), w(t), h[t], v[t]\right)\right) G_2(r)$$
$$\bar{\Omega} := 2(1+\gamma)(1+G_2(r))\exp(\beta r) \tag{4.72}$$

where $G_i : [0, R) \to (0, +\infty)$ $(i = 1, 2)$ are the non-decreasing functions involved in (3.43). Estimate (3.29) with $\bar{M} = \sqrt{\bar{\Omega}}$ and $\bar{\lambda} = \frac{\bar{\omega}}{2}$ is a consequence of estimate (4.71). Since $G_i : [0, R) \to (0, +\infty)$ $(i = 1, 2)$ are non-decreasing functions, we obtain from (4.58), (4.70), (4.72) for all $t \geq 0$:

$$\left\|(\xi(t), w(t), h[t] - h^* \chi_{[0,L]}, v[t])\right\|_X^2 \leq G_1(r) G_2(r) \exp\left(-\frac{\omega t}{\Lambda(r)}\right) \left\|(\xi(0), w(0), h[0] - h^* \chi_{[0,L]}, v[0])\right\|_X^2 \tag{4.73}$$

Estimate (3.28) with $M = \sqrt{G_1(r)G_2(r)}$ and $\lambda = \frac{\omega}{2\Lambda(r)}$ is a direct consequence of (4.73). The proof is complete. ◁

**Proof of Proposition 2:** Let arbitrary $(\xi, w, h, v) \in S$ that satisfies $\left\|(0, w, h - h^* \chi_{[0,L]}, v)\right\|_X \leq \varepsilon$ with $0 < \varepsilon < \min\left(h^*, H_{\max} - h^*\right)/\sqrt{L}$ be given. Using definitions (3.2), (3.3), (3.4), the inequalities $(h(x)v(x) + \mu h'(x))^2 \leq 2h^2(x)v^2(x) + 2\mu^2 (h'(x))^2$, $(w + k\xi)^2 \leq 2w^2 + 2k^2\xi^2$, we obtain:

$$V(\xi, w, h, v) \leq \frac{\delta+2}{2} \int_0^L h(x) v^2(x) dx + \mu^2 \int_0^L h^{-1}(x) (h'(x))^2 dx + \frac{\delta+1}{2} g \left\|h - h^* \chi_{[0,L]}\right\|_2^2 + \frac{3qk^2}{2} \xi^2 + qw^2 \tag{4.74}$$

Since $\int_0^L h(x) dx = m$ (recall definition (3.1)) and since $h^* = m/L$, it follows from continuity of $h$ and the mean value theorem that there exists $x^* \in [0, L]$ such that $h(x^*) = h^*$. Using the Cauchy-Schwarz inequality, we get for all $x \in [0, L]$:

$$h(x) - h^* = \int_{x^*}^x h'(s) ds \Rightarrow \left|h(x) - h^*\right| \leq \int_{\min(x,x^*)}^{\max(x,x^*)} |h'(s)| ds \leq \|h'\|_1 \leq \sqrt{L} \|h'\|_2 \tag{4.75}$$

Consequently, we obtain for all $x \in [0, L]$ from (4.75), (3.13) and the facts that $0 < \varepsilon < \min(h^*, H_{max} - h^*)/\sqrt{L}$, $\|(0, w, h - h^* \chi_{[0,L]}, v)\|_X \leq \varepsilon$:

$$0 < h^* - \varepsilon\sqrt{L} < h^* - \sqrt{L}\|h'\|_2 \leq h(x) \leq h^* + \sqrt{L}\|h'\|_2 \leq h^* + \varepsilon\sqrt{L} < H_{max} \qquad (4.76)$$

Therefore, we obtain from (4.74) and (4.76):

$$\begin{aligned}V(\xi, w, h, v) \leq &\frac{\delta + 2}{2} H_{max} \|v\|_2^2 + qw^2 + \frac{3qk^2}{2}\xi^2 \\ &+ \mu^2 \left(h^* - \varepsilon\sqrt{L}\right)^{-1} \|h'\|_2^2 + \frac{\delta + 1}{2} g \|h - h^* \chi_{[0,L]}\|_2^2\end{aligned} \qquad (4.77)$$

Inequality (3.22) is a direct consequence of (4.77) and definition (3.13).
The proof is complete. ◁

**Proof of Theorem 2:** Let $\sigma, q, k, \delta > 0$, $\omega_1 \in (0, h^*)$, $\omega_2 > 0$ be given constants for which inequalities (3.31), (3.32), (3.35), (3.36) hold. Let $r \in [0, R)$ be a given constant for which inequalities (3.34) hold.

The proof of Lemma 3 shows that the non-decreasing function $\Lambda : [0, R) \to (0, +\infty)$ involved in (3.42) may be defined by (4.25). It follows from (4.25), (3.34) and the fact that $p_2(r) < H_{max}$ (a consequence of the fact that $r \in [0, R)$ and definitions (3.9), (3.11)) that the following inequality holds:

$$\Lambda(r) \leq \frac{1}{2} \max\left(\frac{L^2(\delta + 2)H_{max}}{\pi^2 \omega_1}, (\delta + 1)gL^2 + \frac{2\mu^2}{\omega_1}, qk^2, q\right) \qquad (4.78)$$

It follows from (3.26), (3.34) and definition (3.33) that the following inequality holds:

$$\theta(r) > \tilde{\theta} \qquad (4.79)$$

Due to (4.78), (4.79), (3.34), definition (3.39) and the fact that $p_2(r) < H_{max}$ (a consequence of the fact that $r \in [0, R)$ and definitions (3.9), (3.11)), the functions $\varphi(s)$, $\alpha(s)$ defined by (3.47), (3.48) satisfy the following inequalities:

$$\varphi(r) \geq \omega_1 H_{max} \qquad (4.80)$$

$$\alpha(r) \geq \tilde{\alpha} \qquad (4.81)$$

Inequalities (3.35), (3.36), (3.34), (4.80), (4.81) imply that the following inequalities hold:

$$\begin{aligned}\alpha(r) &> \frac{5\left(H_{max}\tilde{K}^2 + \varepsilon_1\right)}{\delta \mu \gamma} \\ \alpha(r)\frac{\beta}{2} + \frac{\mu\delta\beta\gamma}{4H_{max}}\varphi(r) &> \frac{\varepsilon_2}{p_1^2(r)} \\ \frac{\mu\delta\beta}{8H_{max}}\varphi(r) &> \frac{5L}{6\mu}\end{aligned} \qquad (4.82)$$

Since the functions $p_1(s)$, $\varphi(s)$, $\alpha(s)$ defined by (3.9), (3.47), (3.48), respectively, are continuous functions, inequalities (3.34) and (4.82) imply the existence of $\bar{r} \in (r, R)$ for which the following inequalities hold:

$$p_1(\bar{r}) \geq \omega_1, \quad \sqrt{\frac{2L}{3}\left(\bar{r} + \gamma\bar{r}\exp(\beta\bar{r})\right)} \leq \omega_2 \tag{4.83}$$

$$\alpha(\bar{r}) \geq \frac{5\left(H_{\max}\tilde{K}^2 + \varepsilon_1\right)}{\delta\mu\gamma}$$

$$\alpha(\bar{r})\frac{\beta}{2} + \frac{\mu\delta\beta\gamma}{4H_{\max}}\varphi(\bar{r}) \geq \frac{\varepsilon_2}{p_1^2(\bar{r})} \tag{4.84}$$

$$\frac{\mu\delta\beta}{8H_{\max}}\varphi(\bar{r}) \geq \frac{5L}{6\mu}$$

Consider an arbitrary classical solution of the PDE-ODE system (2.7)-(2.10), (2.12) with (3.27) with $v \in C^0\left(\Re_+; H^1(0,L)\right) \cap C^1\left((0,+\infty); H^1(0,L)\right)$ that satisfies $U(\xi(0), w(0), h[0], v[0]) \leq r + \gamma r \exp(\beta r)$. We make the following claim.

<u>Claim:</u> If $U(\xi(t), w(t), h[t], v[t]) \leq \bar{r} + \gamma\bar{r}\exp(\beta\bar{r})$ for some $t > 0$ then $\frac{d}{dt}U(\xi(t), w(t), h[t], v[t]) \leq 0$.

<u>Proof of Claim:</u> Definition (3.5) and the fact that $U(\xi(t), w(t), h[t], v[t]) \leq \bar{r} + \gamma\bar{r}\exp(\beta\bar{r})$ implies that

$$V(\xi(t), w(t), h[t], v[t]) \leq \bar{r} \tag{4.85}$$

$$\|v_x[t]\|_2^2 \leq 2\bar{r} + 2\gamma\bar{r}\exp(\beta\bar{r}) \tag{4.86}$$

It follows from (4.41) that

$$\|v[t]\|_\infty \leq \sqrt{\frac{2L}{3}\left(\bar{r} + \gamma\bar{r}\exp(\beta\bar{r})\right)} \tag{4.87}$$

Inequalities (3.12), (4.83), (4.85), (4.87) and definitions (3.30), (3.46) imply that

$$\bar{K}(t) \leq \tilde{K} \tag{4.88}$$

The fact that the functions $\varphi(s)$, $\alpha(s)$ defined by (3.9), (3.47), (3.48), are non-increasing functions, in conjunction with (3.31), (3.49), (4.84), (4.85), (4.88), imply that $\frac{d}{dt}U(\xi(t), w(t), h[t], v[t]) \leq 0$. The proof of the claim is complete.

Continuity of the mapping $t \to U(\xi(t), w(t), h[t], v[t])$ (recall that $v \in C^0\left(\Re_+; H^1(0,L)\right)$, $h \in C^1\left(\Re_+ \times [0,L]; (0,+\infty)\right)$ and $v \in C^0\left(\Re_+ \times [0,L]\right)$), the facts that $\bar{r} > r$, $U(\xi(0), w(0), h[0], v[0]) \leq r + \gamma r \exp(\beta r) < \bar{r} + \gamma\bar{r}\exp(\beta\bar{r})$ and the Claim allow us to use a standard contradiction argument and prove that the estimate $U(\xi(t),w(t),h[t],v[t]) \leq \bar{r} + \gamma\bar{r}\exp(\beta\bar{r})$ holds for all $t \geq 0$. Therefore, the Claim in conjunction with continuity of the mapping $t \to U(\xi(t),w(t),h[t],v[t])$ and the fact that $U(\xi(0),w(0),h[0],v[0]) \leq r + \gamma r\exp(\beta r)$ implies that the following estimate holds for all $t \geq 0$:

$$U(\xi(t),w(t),h[t],v[t]) \leq U(\xi(0),w(0),h[0],v[0]) \leq r + \gamma r\exp(\beta r) \tag{4.89}$$

Inequality (4.89) in conjunction with (3.17) shows that $(\xi(t),w(t),h[t],v[t]) \in X_U(r)$ for all $t \geq 0$. Working exactly as in the proof of the claim above and exploiting (4.89), we conclude that (4.88) holds for all $t \geq 0$ as well as the following inequality for all $t \geq 0$:

$$V(\xi(t),w(t),h[t],v[t]) \leq r \tag{4.90}$$

Using (3.44), (3.49), (4.82), (4.88), (4.90), (3.31), the fact that $p_2(r) < H_{\max}$, the facts that $p_1(s)$, $\varphi(s)$, $\alpha(s)$ are non-increasing and $p_2(s)$ is non-decreasing, we get the following differential inequalities for all $t \geq 0$:

$$\begin{aligned}\frac{d}{dt}V(\xi(t),w(t),h[t],v[t]) &\leq -\frac{\mu g}{4}\|h_x[t]\|_2^2 - qk^3\xi^2(t) \\ &\quad -q(q\theta(r)-k)(w(t)+k\xi(t))^2 - \frac{\mu\delta}{2}\int_0^L h(t,x)v_x^2(t,x)dx\end{aligned} \tag{4.91}$$

$$\begin{aligned}\frac{d}{dt}U(\xi(t),w(t),h[t],v[t]) &\leq -\alpha(r)V(t) \\ &\quad -\left(\frac{\delta\gamma}{H_{\max}}\varphi(r)+\frac{\pi^2}{L^2}\right)\frac{\mu}{4}\exp(\beta V(t))\|v_x[t]\|_2^2 \\ &\quad -\left(\alpha(r)-\frac{5(H_{\max}\tilde{K}^2+\varepsilon_1)}{\delta\mu\gamma}\right)\gamma\exp(\beta V(t))V(t)\end{aligned} \tag{4.92}$$

It follows from (3.42), (3.5), (4.90), (4.91), (4.92), we get the following differential inequalities for all $t \geq 0$:

$$\frac{d}{dt}V(\xi(t),w(t),h[t],v[t]) \leq -\omega V(\xi(t),w(t),h[t],v[t]) \tag{4.93}$$

$$\frac{d}{dt}U(\xi(t),w(t),h[t],v[t]) \leq -\bar{\omega}U(\xi(t),w(t),h[t],v[t]) \tag{4.94}$$

where

$$\omega = \frac{1}{\Lambda(r)}\min\left(\frac{\mu g}{4},qk^3,q(q\theta(r)-k),\frac{\mu\delta}{2}\right),$$

$$\bar{\omega} = \min\left(\alpha(r)-\frac{5(H_{\max}\tilde{K}^2+\varepsilon_1)}{\delta\mu\gamma},\left(\frac{\delta\gamma}{H_{\max}}\varphi(r)+\frac{\pi^2}{L^2}\right)\frac{\mu}{2}\right)$$

From this point the proof becomes identical with the proof of Theorem 1. ◁

# 5. Concluding Remarks

This paper showed that a simple modification of the nonlinear feedback law proposed in [22], namely an increase of the gain of the liquid momentum, can guarantee robustness with respect to wall friction. The constructed stabilizing feedback law does not require exact knowledge of the friction coefficient (for which there is considerable uncertainty). The obtained robustness results are not a surprise: Lyapunov function(al)s are among the most important tools of guaranteeing robustness and the feedback laws were designed by means of the CLF methodology. Therefore, it is justified to claim that robustness is an intrinsic property for the CLF methodology-a well-known fact in the literature of finite-dimensional systems. The results of the present paper confirm this fact for the infinite-dimensional case of systems described by PDEs.

The present paper did not study the existence/uniqueness issue for the solutions of the closed-loop system. This is a completely different topic for future research. To this purpose, ideas contained in [24,31] can be used. However, the obtained stability estimates will certainly help the analysis. Another open issue is the construction of numerical schemes for the numerical approximation of the solutions of the closed-loop system. The issue was studied in [23], where we noticed the lack of numerical schemes that can be applied to the closed-loop system. Again, the stability estimates provided in the present work can help the analysis. Even in the preliminary numerical study contained in [23] the time step was selected based on the discretized version of the CLF. We expect that the Lyapunov functionals that were provided in the present work will help the stabilization of the numerical schemes that will be used for the approximation of the solutions of the closed-loop system.

**References**

[1] Barré de Saint-Venant, A. J. C., "Théorie du Mouvement non Permanent des Eaux, avec Application aux Crues des Rivières et a l'Introduction de Marées dans Leurs Lits", *Comptes Rendus de l'Académie des Sciences*, 73, 1871, 147–154 and 237–240.

[2] Bastin, G., J.-M. Coron and B. d'Andréa Novel, "On Lyapunov Stability of Linearised Saint-Venant Equations for a Sloping Channel", *Networks and Heterogeneous Media*, 4, 2009, 177-187.

[3] Bastin, G. and J.-M. Coron, *Stability and Boundary Stabilization of 1-D Hyperbolic Systems*, 2016, Birkhäuser.

[4] Bastin, G., J.-M. Coron and A. Hayat, "Feedforward Boundary Control of 2x2 Nonlinear Hyperbolic Systems with Application to Saint-Venant Equations", *European Journal of Control*, 57, 2021, 41-53.

[5] Bresch, D. and B. Desjardins, ''On the Construction of Approximate Solutions for the 2D Viscous Shallow Water Model and for Compressible Navier–Stokes Models'', *Journal de Mathématiques Pures et Appliquées*, 86, 2006, 362-368.

[6] Bresch, D. and P. Noble, "Mathematical Justification of a Shallow Water Model", *Methods and Applications of Analysis*, 14, 2007, 87–118.

[7] Chin-Bing S. A., P. M. Jordan and A. Warn-Varnas, "A Note on the Viscous, 1D Shallow Water Equations: Traveling Wave Phenomena", *Mechanics Research Communications*, 38, 2011, 382-387.

[8] Coron, J. M., "Local Controllability of a 1-D Tank Containing a Fluid Modeled by the Shallow Water Equations", *ESAIM: Control Optimisation and Calculus of Variations*, 8, 2002, 513–554.

[9] Coron, J.-M., *Control and Nonlinearity*, Mathematical Surveys and Monographs, Volume 136, American Mathematical Society, 2007.

[10] Coron, J.-M., B. d'Andréa Novel and G. Bastin, "A Strict Lyapunov Function for Boundary Control of Hyperbolic Systems of Conservation Laws", *IEEE Transactions on Automatic Control*, 52, 2007, 2–11.

[11] Coron, J.-M., A. Hayat, S. Xiang and C. Zhang. "Stabilization of the Linearized Water Tank System", *Archive for Rational Mechanics and Analysis*, 244, 2022, 1019–1097.

[12] de Halleux, J. and G. Bastin, "Stabilization of Saint-Venant Equations Using Riemann Invariants: Application to Waterways with Mobile Spillways", *IFAC Proceedings*, 35, 2002, 131-136.

[13] Diagne, A., M. Diagne, S. Tang and M. Krstic, "Backstepping Stabilization of the Linearized Saint-Venant–Exner Model", *Automatica*, 76, 2017, 345–354.

[14] Diagne, M., S.-X. Tang, A. Diagne and M. Krstic, "Control of Shallow Waves of two Unmixed Fluids by Backstepping", *Annual Reviews in Control*, 44, 2017, 211–225.

[15] Dos Santos, V., G. Bastin, J.-M. Coron and B. d'Andréa-Novel, "Boundary Control with Integral Action for Hyperbolic Systems of Conservation Laws: Stability and Experiments", *Automatica*, 44, 2008, 1310-1318.

[16] Dubois, F., N. Petit and P. Rouchon, "Motion Planning and Nonlinear Simulations for a Tank Containing a Fluid", *Proceedings of the 1999 European Control Conference (ECC)*, 1999, 3232-3237.

[17] Gerbeau, J.-F. and B. Perthame, "Derivation of Viscous Saint-Venant System for Laminar Shallow-Water: Numerical Validation", *Discrete and Continuous Dynamical Systems, Series B*, 1, 2001, 89–102.

[18] Hayat, A. and P. Shang, "A Quadratic Lyapunov Function for Saint-Venant Equations with Arbitrary Friction and Space-Varying Slope", *Automatica*, 100, 2019, 2–60.

[19] Hayat, A. and P. Shang, "Exponential Stability of Density-Velocity Systems with Boundary Conditions and Source Term for the $H^2$ Norm", *Journal de Mathematiques Pures et Appliquees*, 153, 2021, 187-212.

[20] Karafyllis, I. and M. Krstic, "Global Stabilization of a Class of Nonlinear Reaction-Diffusion PDEs by Boundary Feedback", *SIAM Journal on Control and Optimization*, 57, 2019, 3723-3748.

[21] Karafyllis, I. and M. Krstic, "Global Stabilization of Compressible Flow Between Two Moving Pistons", *SIAM Journal on Control and Optimization*, 60, 2022, 1117-1142.

[22] Karafyllis, I. and M. Krstic, "Spill-Free Transfer and Stabilization of Viscous Liquid", to appear in *IEEE Transactions on Automatic Control* (see also arXiv:2108.11052 [math.OC]).

[23] Karafyllis, I., F. Vokos and M. Krstic, "Output-Feedback Control of Viscous Liquid-Tank System and its Numerical Approximation", submitted to *Automatica* (see also arXiv:2201.13272 [math.OC]).

[24] Kloeden, P. E., "Global Existence of Classical Solutions in the Dissipative Shallow Water Equations", *SIAM Journal on Mathematical Analysis*, 16, 1985, 301-315.

[25] Lannes, D., *The Water Waves Problem. Mathematical Analysis and Asymptotics*, 2013, American Mathematical Society.

[26] Litrico, X. and V. Fromion, "Boundary Control of Linearized Saint-Venant Equations Oscillating Modes", *Automatica*, 42, 2006, 967–972.

[27] Mascia, C. and F. Rousset, "Asymptotic Stability of Steady-States for Saint-Venant Equations with Real Viscosity", in C. Calgaro, J.-F. Coulombel and T. Goudon (Eds), *Analysis and Simulation of Fluid Dynamics*, Advances in Mathematical Fluid Mechanics, 2006, Birkhäuser, 155–162.

[28] Petit, N. and P. Rouchon, "Dynamics and Solutions to Some Control Problems for Water-Tank Systems", *IEEE Transactions on Automatic Control*, 47, 2002, 594-609.

[29] Prieur, C. and J. de Halleux, "Stabilization of a 1-D Tank Containing a Fluid Modeled by the Shallow Water Equations", *Systems & Control Letters*, 52, 2004, 167-178.

[30] Shames, I. H., *Mechanics of Fluids*, 2nd Edition, McGraw-Hill International Editions, 1989.

[31] Sundbye, L., "Global Existence for the Dirichlet Problem for the Viscous Shallow Water Equations", *Journal of Mathematical Analysis and Applications*, 202, 1996, 236-258.

[32] Vazquez, R. and M. Krstic, *Control of Turbulent and Magnetohydrodynamic Channel Flows*, Birkhäuser, 2007.